\documentclass[final]{article}

\usepackage{graphicx, amssymb, amsmath, epsfig}
\usepackage{color}
\usepackage{pst-text}

\graphicspath{{./}{figs/}}

\newcommand{\R}{\mathbb R}

\newcommand{\Gcoarse}{\widetilde G}
\newcommand{\Ncoarse}{\widetilde N}
\newcommand{\Ncoarsebis}{$\Ncoarse \phantom{^{3^3}}$ \hskip-11pt}
\newcommand{\kcoarse}{\widetilde k}

\newcommand{\xcoarse}{\widetilde x}
\newcommand{\ucoarse}{\widetilde U}
\newcommand{\ba}{\begin{array}{cc}}
\newcommand{\ea}{\end{array}}

\newtheorem{dfn}{Definition}[section]
\newtheorem{teo}{Theorem}[section]

\title{A Patchy Dynamic Programming Scheme for a Class of Hamilton-Jacobi-Bellman Equations\thanks{
The authors wish to acknowledge the support obtained by the following grants: AFOSR Grant no.\ FA9550-10-1-0029, 
ITN - Marie Curie Grant n. 264735-SADCO 
and SAPIENZA 2009 {``Analisi ed approssimazione di modelli differenziali nonlineari in fluidodinamica e scienza dei materiali''}. The authors also wish to thank the CASPUR Consortium for its technical support.}
}

\author{
Simone Cacace\footnote{Dipartimento di Matematica, 
SAPIENZA - Universit\`a di Roma, P.le Aldo Moro 2, 00185 Rome, Italy. 
\texttt{simone.cacace@gmail.com}}, 
Emiliano Cristiani\footnote{Dipartimento di Matematica, 
SAPIENZA - Universit\`a di Roma, P.le Aldo Moro 2, 00185 Rome, Italy. 
\texttt{emiliano.cristiani@gmail.com}},
Maurizio Falcone\footnote{(corresponding author) Dipartimento di Matematica, 
SAPIENZA - Universit\`a di Roma, P.le Aldo Moro 2, 00185 Rome, Italy. 
\texttt{falcone@mat.uniroma1.it}}, 
Athena Picarelli\footnote{Dipartimento di Matematica, 
SAPIENZA - Universit\`a di Roma, P.le Aldo Moro 2, 00185 Rome, Italy. 
\texttt{athena.picarelli@gmail.com}} 
}

\begin{document}

\maketitle

\begin{abstract}
In this paper we present a new algorithm for the solution of Hamilton-Jacobi-Bellman equations related to optimal control problems. The key idea is to divide the domain of computation into subdomains which are shaped by the optimal dynamics of the underlying control problem. This can result in a rather complex geometrical subdivision, but has the advantage that every subdomain is invariant with respect to the optimal dynamics and then the solution can be computed \textit{independently} in each subdomain. The features of this dynamics-dependent domain decomposition can be exploited to speed up the computation and for an efficient parallelization, since the classical transmission conditions at the boundaries of the subdomains can be avoided. For their properties, the subdomains are \textit{patches} in the sense introduced by Ancona and Bressan in 1999. Several examples in dimension two and three illustrate the properties of the new method.
\end{abstract}


\textbf{AMS.}
65N55, 49L20

\pagestyle{myheadings}
\thispagestyle{plain}
\markboth{S. CACACE, E. CRISTIANI, M. FALCONE, A. PICARELLI}
{A Patchy Dynamic Programming Scheme for HJB Equations}


\section{Introduction}
The numerical solution of partial differential equations (PDEs) obtained by applying the dynamic programming principle (DPP) to nonlinear optimal control problems is a challenging topic that can have a great impact 
in many areas, e.g.\ robotics, aeronautics, electrical and aerospace engineering. 
Indeed, by means of the DPP one can characterize the value function of a fully--nonlinear control problem (including also state/control constraints) as the unique viscosity solution of a nonlinear Hamilton--Jacobi 
equation, and, even more important, from the solution of this equation one can derive the approximation of an optimal feedback control. 
This result is the main motivation for the PDE approach to control problems and represents the main advantage over other methods, such as those based on the Pontryagin minimum principle. 
It is worth mentioning that the characterization via the Pontryagin principle gives only necessary conditions for the optimal trajectory and optimal open-loop control. In addition, the numerical procedures for solving the associated system of ordinary differential equations can be very complicated. In real applications, a good initial guess for the co-state often requires a long and tedious trial-error procedure to be found. 
This is why it can be useful to combine the DPP and the Pontryagin  approaches, using the approximate value function to compute a 
suitable initial guess for the co-state, as proposed in \cite{CM10}. 

In this paper we mainly focus on the \emph{minimum time problem}, which is associated to the following Hamilton-Jacobi-Bellman equation
\begin{equation}\label{HJB}
\left\{
\begin{array}{ll}
\displaystyle\max_{a\in A}\left\{-f(x,a)\cdot\nabla u(x)-1\right\}=0\,, & x\in\R^d\backslash\Omega_0 \\
u(x)=0\,, & x\in\Omega_0
\end{array}
\right.
\end{equation}
where $d$ is the dimension of the state, $A\subset\R^m$ is the set of admissible controls, $\Omega_0$ is the target to be reached in minimal time and $f:\R^d\times A\to\R^d$ is the dynamics of the system. 
The value function $u:\R^d\rightarrow \R$ at the point $x$ is the minimal time to reach the target starting from $x$ (note that $u(x)=+\infty$ if the target is not reachable).
For numerical purposes, the equation is solved in a bounded domain $\Omega\supset\Omega_0$, so that also boundary conditions on $\partial\Omega$ are needed. A rather standard choice when one has not  additional information on the solution and deals with target problems is imposing state constraint boundary conditions. 

The techniques used to obtain a numerical approximation of the viscosity solution of equation (\ref{HJB}) have been mainly based on Finite Differences \cite{CL84, S99} and Semi-Lagrangian schemes \cite{F97,  FFbook}. 
More recently, Finite Elements methods based on Discontinuous Galerkin approximations have been proposed, due to their ability to deal with non regular functions, 
which is the typical case in the framework of viscosity solutions \cite{CC07, CS07, S07}. 
It is important noting that traditional approximation schemes presented, for example, in \cite{CL84} and \cite{F97}, are based on fixed-point iterations, meaning that the solution is computed at each node of the grid at every iteration until convergence. 
Denoting by $M$ the number of nodes in each dimension and considering that the number of iterations needed for convergence is of order $O(M)$, the total cost of these full-grid schemes is $O(M^{d+1})$. 
We easily conclude that classical algorithms are very expensive when the state dimension is $d\geq 3$, although they are rather efficient for low-dimensional control problems as shown 
in \cite{F97} (see also the book \cite{FFbook}). 

The ``curse of dimensionality'' has been attacked in many ways and new techniques have been proposed to accelerate convergence and/or to reduce memory allocation. 
In \cite{BCZ10} authors proposed an algorithm that allows one to allocate only a small portion of the grid at every iteration. 
Another proposal to reduce the computational effort is given by the Fast Marching method, introduced in \cite{S96,T95} for the Eikonal equation. While the full-size grid is always allocated, the computation is restricted to a small portion of the grid, thus saving CPU time. The cost of this method is of order $O(M^d\log M^d)$.
Despite the efficiency of the Fast Marching method, at present its application to more general equations of the form \eqref{HJB} is not an easy task and it is still under investigation \cite{CCF11,CFFM08,C09,SV03}. 

Other methods have been proposed exploiting the idea that one can accelerate convergence by alternating the order in which the grid nodes are visited giving rise to the so-called ``sweeping methods''. Unlike Fast Marching methods, these methods do not require a special ordering of the grid nodes and are somehow blind, so it could be difficult to prove that they converge after a finite number of sweeps. However, they are easy to implement and they have been shown to be efficient for the Eikonal equation \cite{Z05} and, more recently, for rather general Hamiltonians \cite{TCOZ04}.

Another strategy is based on the decomposition of the domain $\Omega$. The problem is actually solved in subdomains $\Omega^j$, $j=1,\ldots,R$, 
whose size is chosen in order to reduce the number of grid nodes to a manageable size.
Therefore, rather than solving a unique huge problem, one can solve $R$ smaller subproblems working simultaneously on several processors.  
Depending on the choice of the subdomains $\Omega^j$, we can have overlapping regions or interfaces between the subdomains. This is a delicate point, since at each iteration it is necessary exchanging information between processors to make subdomains communicate with each other. Without this communication the result will not be correct. The interested reader can find in the book \cite{QV99} a comprehensive introduction to domain decomposition techniques, whereas for an application to Hamilton-Jacobi equations we refer to \cite{CFLS94, FLS94}.

Finally, a decomposition of the domain based on the concept of ``patchy feedbacks'' has been proposed. It was introduced by Ancona and Bressan in \cite{AB99}, where the authors studied 
the problem of the asymptotic stabilization of a control system. Their main result (see Theorem \ref{ancona_bressan} in Section \ref{patchy_feedback}) states that, under suitable assumptions 
on the control system, stabilization can be obtained by means of a special feedback control which is piecewise constant on a particular partition of the domain. 
Such a partition has the fundamental property that each part, or ``patch'', is positive-invariant with respect to the optimal dynamics driving the system. This is the spirit of what we call ``patchy method''. Unfortunately, the result of Ancona and Bressan is purely theoretical and their patchy decomposition turns out to be not constructive. 
Thus, one has to face the problem of a numerical approximation of such a dynamics-invariant domain decomposition.

A first example of discrete patchy method has been proposed by Navasca and Krener in \cite{KN07, KN08}. The authors adopt a formal method developed by Al'brekht \cite{A61}  that essentially translates the Hamilton-Jacobi-Bellman equation associated to a control problem into a system of algebraic equations, whose unknowns are the coefficients of the expansions in power series of the cost and  optimal feedback. This gives an approximate solution in a small neighbourhood of the origin, which is the first patch of their domain decomposition. 
The solution is then extended to new patches around the first one, by picking some boundary points from where optimal trajectories emanate (they are computed numerically backward in time). 
Those points define the centres of new neighbourhoods that can be used to restart the method. 
The solution is then obtained iteratively by fitting together the approximations in all the patches. More recently, it has been shown that this technique can be extended to obtain high-order accuracy in the regions where 
the value function is smooth (see \cite{H11} for more details).

Despite the high speed of the method, that actually does not use any grid, there are many open questions on its application. The first limitation is the strong regularity assumptions on the solution necessary to set the problem in these terms. Indeed, it is well known that the most simple control problems may have optimal controls which are not even continuous \cite{BCDbook}. The second crucial point is the construction of the patchy decomposition that, in the examples contained in \cite{KN07, KN08}, appears not to be completely invariant with respect to the optimal dynamics. This makes the solution rather inaccurate, especially near the boundaries of the patches.

The goal of this paper is to present and investigate a new patchy technique based on a semi-Lagrangian scheme that leads to a dynamics-dependent partition of the domain $\Omega$. Subdomains turn out to be, up to a discretization error, \textit{invariant with respect to the optimal dynamics}, meaning that optimal dynamics do not cross the boundaries of the subdomains or, equivalently, that boundaries of the subdomains are optimal trajectories to the target. 
The algorithm consists of three main steps: first, a rough solution of \eqref{HJB} is computed on a coarse grid. Then, the feedback optimal control is used to obtain the dynamics-dependent domain decomposition. Finally, the solution of the equation is computed on a (much) finer grid, \textit{independently} in each subdomain, by serial or parallel computation. We will see that the invariance of the subdomains can be exploited to lower the iterations needed to reach convergence and to parallelize efficiently the computation on distributed-memory architectures, since there is no need of communication among processors.

The paper is organized as follows. 
Section 2 is devoted to the presentation of the semi-Lagrangian scheme and the classical domain decomposition technique for equation \eqref{HJB}. Moreover, we briefly describe the results on patchy methods which have been proved by Ancona and Bressan \cite{AB99}. 

In section 3 we present the patchy domain decomposition
method and our algorithm to split the domain into invariant subdomains. We discuss
there several issues related to the implementation of the method and its parallelization.
Section 4 is devoted to the numerical tests on control problems in dimension two, and
section 5 presents some improvements of the basic algorithm. Finally, in section 6 we
present some tests in dimension three.
%
%
%
%
%
%
%
%
%
%
%
%
\section{Background}
In this section we briefly introduce the numerical scheme used to discretize equation \eqref{HJB} and the classical domain decomposition technique, including an algorithm that will be used in the following. 
Next, we recall the notion of patchy decomposition and the result by Ancona and Bressan \cite{AB99} concerning the asymptotic stabilization of control systems by means of patchy feedbacks, that inspired our patchy numerical method. 
\subsection{The semi-Lagrangian scheme}
We introduce a structured grid $G$ on $\Omega$ with nodes $x_i,\ i=1,\ldots,N$. We also denote by $\mathring{G}$ the internal nodes of $G$ 
and by $\partial G$ its boundary, whose nodes act as \emph{ghost nodes}.
We map all the values at the nodes onto a $N$-dimensional vector $U=(U_1,\ldots, U_{N})$. 
Let us denote by $h_{i,a}>0$ a (fictitious) time step, possibly depending on the node $x_i$ and control $a$, (see the book \cite{FFbook} for details) and by $k>0$ the space step. 
By a standard semi-Lagrangian discretization \cite{BF90, BF90b, F97} of  \eqref{HJB}, 
 it is possible to obtain the following scheme in fixed-point form 
\begin{equation}\label{SLscheme}
U=F(U)\,,
\end{equation}
where $F:\R^{N}\rightarrow\R^{N}$ is defined componentwise  by 
$$
\left[F(U)\right]_i=\left\{
\begin{array}{ll}
\displaystyle\min_{a\in A}\left\{I\left[U\right](x_i+h_{i,a}f(x_i, a))+h_{i,a}\right\} & x_i\in\mathring{G}\setminus\Omega_0\,,\\
0 & x_i\in\Omega_0\cap G\,,\\
+\infty & x_i \in\partial G\,.
\end{array}
\right.
$$
The discrete value function $U$ is extended on the whole space $\Omega$ by the interpolation operator $I$. In order to fix ideas, one should think to the linear interpolation in $\R^d$ described in \cite{CFF04} but other choices are available \cite{FFbook}.

We choose the time step $h_{i,a}$ such that $|h_{i,a} f(x_i, a)|=k$ for every $i=1,\ldots, N$ and $a\in A$, so that the point $x_i+h_{i,a}f(x_i, a)$ falls in one of the first neighboring cells around $x_i$. 
The minimum over $A$ is evaluated by direct comparison, discretizing the set $A$ with $N_c$ points. Note that defining $F(U)=+\infty$ on $\partial G$ corresponds to imposing state constraint boundary conditions. 
The final iterative scheme reads
\begin{equation}\label{iterativeSLscheme}
 U^{(n+1)}=F(U^{(n)})\,,\qquad U^{(0)}=\left\{\begin{array}{ll}
                                0 & \mbox{on }\Omega_0\cap G \\
                                +\infty & \mbox{otherwise}
                               \end{array}\right.
.
\end{equation}
We refer to  \cite{F97, FFbook} for details and convergence results.
With the discrete value function $U$ in hand, we can obtain a discrete feedback map $a^*_k:\Omega\rightarrow A$ just defining
\begin{equation}\label{feedh}
a^*_k(x):= \hbox{arg}\min\limits_{a\in A} \{I[U](x+h_{x,a}f(x, a))+h_{x,a}\}\,. 
\end{equation}
Under rather general assumptions \cite{F01}, it can be shown that this is an approximation of the feedback map constructed for the continuous problem as
\begin{equation} \label{feed}
a^*(x):= \hbox{arg}\max\limits_{a\in A} \{- f(x,a)\cdot  \nabla u(x)-1\}\,. 
\end{equation}
A detailed discussion on the construction of the feedback maps and of the optimal trajectories solution of the minimum time problem can be found in \cite{F97,F01}. 
%
%
%
%
%
%
\subsection{Domain Decomposition method}
\label{sec:dd_algorithm}
The domain decomposition method allows one to split the problem in $\Omega$ into $R$ subproblems in subsets
$\Omega^j,\ j=1,\ldots,R$ such that $\Omega=\Omega^1\cup\ldots\cup\Omega^R$. 
For every pair $j\neq \ell$ of indices corresponding to adjacent subdomains $\Omega^j$ and $\Omega^\ell$, 
let us denote by $\Omega^{j\ell}$  the non-empty overlapping zone $\Omega^j\cap\Omega^\ell$, which is assumed to contain at least one grid cell. 

We also denote by $N^j$ the number of nodes of $\Omega^j$, by $U^j$ the restriction of $U$ to $\Omega^j$ 
and by $F^j:\R^{N^j}\to\R^{N^j}$ the restriction of the operator $F$ in \eqref{SLscheme} to $\Omega^j$. Then, we define 
globally in $\Omega$ the following splitting operator $F_{\text{\tiny{SPLIT}}}:\R^N\to\R^N$, given componentwise by 
$$
[F_{\text{\tiny{SPLIT}}}(U)]_i\equiv
\left\{
\begin{array}{ll}
[F^j(U^j)]_i & \mbox{if } \exists j\mbox{ such that } x_i\in \Omega^j\setminus \bigcup_{\ell\neq j} \Omega^{j\ell}\,,\\
\displaystyle\min_{j\,:\,x_i\in\Omega^j}\{[F^j(U^j)]_i\} & \mbox{otherwise}\,.
\end{array}
\right.
$$
Following \cite{FLS94}, it is easy to prove that fixed-point iterations for $F$ and $F_{\text{\tiny{SPLIT}}}$ lead to the same solution.

We now describe a simple algorithm to compute the fixed point of $F_{\text{\tiny{SPLIT}}}$.\\

\noindent{\bf Domain Decomposition Algorithm:}
\begin{enumerate}
\item[Step 1.](Initialization)  For $n=0$ the initial guess $U^{(0)}\in \mathbb{R}^N$ is fixed to $0$ on the nodes corresponding to the target $\Omega_0$ and $+\infty$ elsewhere.\\
\item[Step 2.](Computation) $U^{(n+{1/2})}$ is computed separately in every subdomain $\Omega^j$ by 
$$U^{j,(n+{1/2})}=F^j(U^{j,(n)})\qquad j=1,\ldots, R\,.$$
\item[Step 3.](Coupling) 
$$
U^{(n+1)}_i=
\left\{ 
\begin{array}{ll}
U^{j,(n+{1/2})}_i & \mbox{if } \exists j\mbox{ such that } x_i\in \Omega^j\setminus \bigcup_{\ell\neq j} \Omega^{j\ell}\,,\\
\displaystyle\min_{j\,:\,x_i\in\Omega^j}\left\{U^{j,(n+{1/2})}_i\right\} & \mbox{otherwise}\,. 
\end{array}
\right.
$$
\item[Step 4.] (Stopping criterion) If $\|U^{(n+1)}-U^{(n)}\|_\infty>$\emph{tol} go to Step 2 with $n\leftarrow n+1$, otherwise stop.\\ 
\end{enumerate}

In order to speed up the convergence of the algorithm above, we use iterations of Gauss-Seidel type, meaning that we employ the updated values of the nodes as soon as they are available. 

The domain decomposition method can be used in a natural way to parallelize the computation. Simply, each subdomain is assigned to  a processor. Synchronization among processors is performed at each iteration after step 3. Note that, as a limit case, the algorithm can be also used in serial computation (one processor).  

From now on the final solution computed by the domain decomposition algorithm will be denoted by $U_{\!D\!D}$.
\subsection{Patchy feedbacks and stabilization for control systems} \label{patchy_feedback} 
The inspiring idea at the basis of the numerical method we present in the next section is the notion of a ``patch''. 
It has been introduced by Ancona and Bressan in \cite{AB99} in the context of the stabilization of control systems. 
Let us recall here for completeness their main definitions and results.
We refer the interested reader to \cite{AB99} for details and to \cite{AB03, AB07, B98, BP08} for further results on this subject.

We consider the control system
\begin{equation}\label{sist_contr}
\dot{y}(t)=f(y(t), a(t)) \qquad a(t)\in A,
\end{equation}
assuming that the control set $A\subset\R^m$ is compact and the dynamics $f:\R^d\times A\rightarrow\R^d$ is sufficiently smooth. 
Moreover we choose as admissible controls all the functions $a$ belonging to
$$
\mathcal A:=\left\{a:(0, \infty)\rightarrow A\ | \  a \mbox{ is measurable}\right\}.
$$
For every point $x\in\R^d$ and admissible control $a_0\in \mathcal A$, we denote by $y(\cdot\, ; x, a_0)$ the absolutely continuous function defined on some maximal interval $[0; \tau^{\max}(x,a_0))$ satisfying the system (\ref{sist_contr}) with initial condition $y(0)=x$ and control $a_0$.

The following definition extends to control systems the classical notion of stability.
\begin{dfn} 
The system (\ref{sist_contr}) is said to be globally {\em asymptotically controllable} (to the origin) if the following holds:
\begin{enumerate}
\item for each $x\in\R^d$ there exists some admissible control $a_0$ such that the trajectory $t\rightarrow y(t)=y(t; x, a_0)$ is defined for all $t\geq 0$ (i.e.\ $\tau^{\max}(x,a_0)=\infty$) and $y(t)\rightarrow 0$ as $t\rightarrow\infty$;
\item for each $\varepsilon>0$ there exists $\delta>0$ such that for each $x\in\R^d$ with $|x|<\delta$ there is an admissible control $a_0$ as in 1. such that $|y(t)|<\varepsilon$ for all $t\geq 0$.
\end{enumerate}
\end{dfn}
\medskip
\noindent
Given an asymptotically controllable system, a classical problem is finding a feedback control $a=\psi(y):\R^d\to A$ such that \emph{all the trajectories} of the corresponding closed loop system
\begin{equation}\label{sist_cl_loop}
\dot{y}=f(y, \psi(y))
\end{equation}
tend asymptotically to the origin. 
Since this problem may not admit any solution in the class of \emph{continuous feedbacks}, Ancona and Bressan introduce and investigate the properties of a particular class of \emph{discontinuous feedbacks}, the so-called\textit{ patchy feedbacks}.

\noindent
The following definition introduces the fundamental concept of a patch.
\begin{dfn}
Let $\mathcal P\subset\R^d$ be an open domain with smooth boundary $\partial\mathcal P$ and $g$ be a smooth vector field defined on a neighborhood of  $\overline{\mathcal P}$. 
We say that the pair $(\mathcal P, g)$ is a {\em patch} if $\mathcal P$ is a positive-invariant region for $g$, i.e.\ at every boundary point $y\in\partial\mathcal P$ the inner product of $g$ with the outer normal $n$ satisfies
$$
\left\langle g(y), n(y)\right\rangle<0.
$$
\end{dfn}

Then, by means of a superposition of patches, we get the notion of a patchy vector field on a domain $\Omega\subset\R^d$.
\begin{dfn}
We say that $g:\Omega\rightarrow\R^d$ is a {\em patchy vector field} if there exists a family of patches $\left\{(\Omega^\alpha, g^{\alpha}): \alpha\in\mathcal I\right\}$ such that
\begin{itemize}
\item $\mathcal I$ is a totally ordered index set,
\item the open sets $\Omega^{\alpha}$ form a locally finite covering of $\Omega$,
\item the vector field $g$ can be written in the form
$$
g(y)=g^{\alpha}(y)\qquad \mbox{if}\qquad y\in\Omega^{\alpha}\setminus\bigcup_{\beta>\alpha}\Omega^{\beta}.
$$
\end{itemize}
\end{dfn}
We use $\left(\Omega, g, (\Omega^\alpha, g^{\alpha})_{\alpha\in\mathcal{I}}\right)$ to denote the patchy vector field and the family of patches.
By applying the previous definitions to the closed loop system (\ref{sist_cl_loop}) we define a \emph{patchy feedback} control as a piecewise constant map $\psi:\R^d\rightarrow A$ 
such that the vector field $g(y):=f(y, \psi(y))$ is a patchy vector field. More precisely:
\begin{dfn}
Let $\left(\Omega, g, (\Omega^\alpha, g^{\alpha})_{\alpha\in\mathcal{I}}\right)$ be a patchy vector field. Assume that there exists control values $\psi^{\alpha}\in A$ such that, for each $\alpha\in\mathcal I$
$$
g^{\alpha}(y)=f(y, \psi^{\alpha})\qquad\forall y\in\Omega^{\alpha}\setminus\bigcup_{\beta>\alpha}\Omega^{\beta}.
$$
Then, the piecewise constant map
$$
\psi(y)=\psi^{\alpha}\qquad\mbox{if}\qquad y\in\Omega^{\alpha}\setminus\bigcup_{\beta>\alpha}\Omega^{\beta}
$$
is called a {\em patchy feedback} control on $\Omega$.
\end{dfn}
\medskip
\begin{dfn}
A patchy feedback control $\psi:\R^d\setminus{\left\{0\right\}}\rightarrow A$ is said to {\em asymptotically stabilize} the closed loop system (\ref{sist_cl_loop}) with respect to the origin if the following holds:
\begin{enumerate}
\item for each $x\in\R^d\setminus{\left\{0\right\}}$ and for every trajectory $y(\cdot)$ of (\ref{sist_cl_loop}) starting from $x$ one has
 $y(t)\rightarrow 0$ as $t\rightarrow\tau^{\max}(x)$,
\item for each $\varepsilon>0$  there exists $\delta>0$ such that, for each $x\in\R^d\setminus{\left\{0\right\}}$ with $|x|<\delta$ and for every trajectory $y(\cdot)$ of (\ref{sist_cl_loop}) starting from $x$ 
one has $|y(t)|<\varepsilon$, for all $0\leq t<\tau^{\max}(x)$.\\
\end{enumerate}
\end{dfn}
\noindent Finally, the main result of Ancona and Bressan can be summarized as follows.
\begin{teo}\label{ancona_bressan}
If the system (\ref{sist_contr}) is asymptotically controllable, then it admits an asymptotically stabilizing patchy feedback control.\\ 
\end{teo}


\section{The patchy domain decomposition}\label{sec:patchy}
In this section we introduce our new numerical method for solving equations of Hamilton-Jacobi-Bellman type. In particular we focus on the minimum time problem (\ref{HJB}). 
The main feature of the new method is the technique we use to construct the subdomains of the decomposition, which 
are (approximate) patches in a sense inspired by the definitions of the previous section. Indeed, we will see that these patches turn out to be (almost) invariant with respect to the optimal dynamics driving the system, meaning that the optimal dynamics do not cross their boundaries. Even if this construction could lead to a rather complicated domain decomposition, it has the clear advantage that we do not need to apply any transmission condition between subdomains.

Let us introduce two rectangular (structured) grids. The first grid is rather \emph{coarse} because it is used for preliminary (and fast) computations only. 
It will be denoted by $\Gcoarse$ and its nodes by $\xcoarse_1,\ldots,\xcoarse_{\Ncoarse}$, where $\Ncoarse$ is the total number of nodes.  
We will denote the space step for this grid by $\kcoarse$ and the approximate solution of the equation (\ref{HJB}) on this grid by $\ucoarse_{P}$. 

The second grid is instead \emph{fine}, being the grid where we actually want to compute the numerical solution of the equation. It will be denoted by $G$ and its nodes by $x_1,\ldots,x_N$, 
where $N$ is the total number of nodes ($N\!\!>>\!\!\Ncoarse$). We will denote the space step for this grid by $k$ and the solution of the equation (\ref{HJB}) on this grid by $U_P$. 
We also choose the number $R$ of subdomains (patches) to be used in the patchy decomposition and we divide the target $\Omega_0$ in $R$ parts denoted by $\Omega_0^j$, with $j=1,\ldots,R$.\\ 
The patchy method can be described as follows.\\

\begin{figure}[h!]
\begin{center}
\begin{tabular}{cc}
\includegraphics[width=0.4\textwidth]{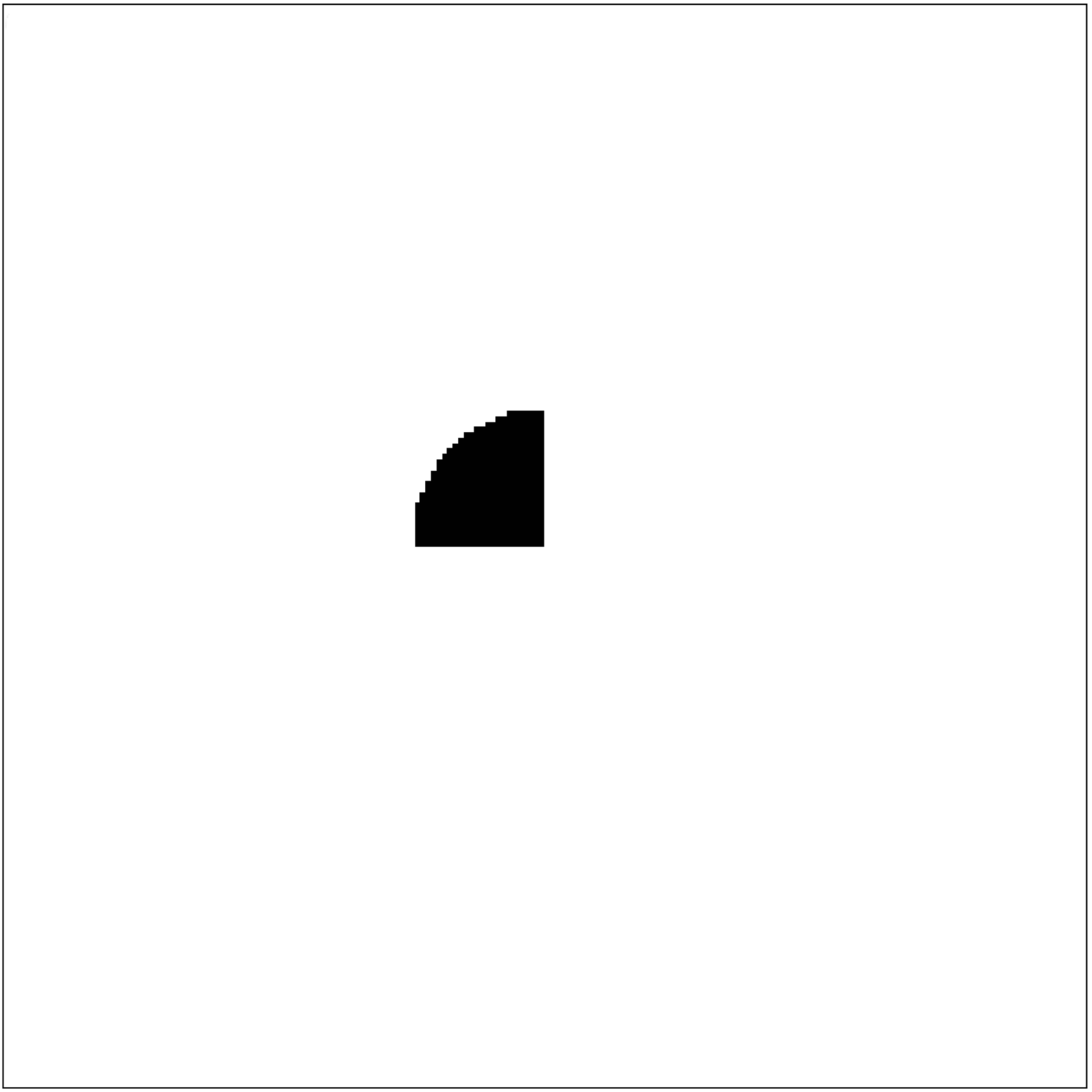} \qquad & \qquad 
\includegraphics[width=0.4\textwidth]{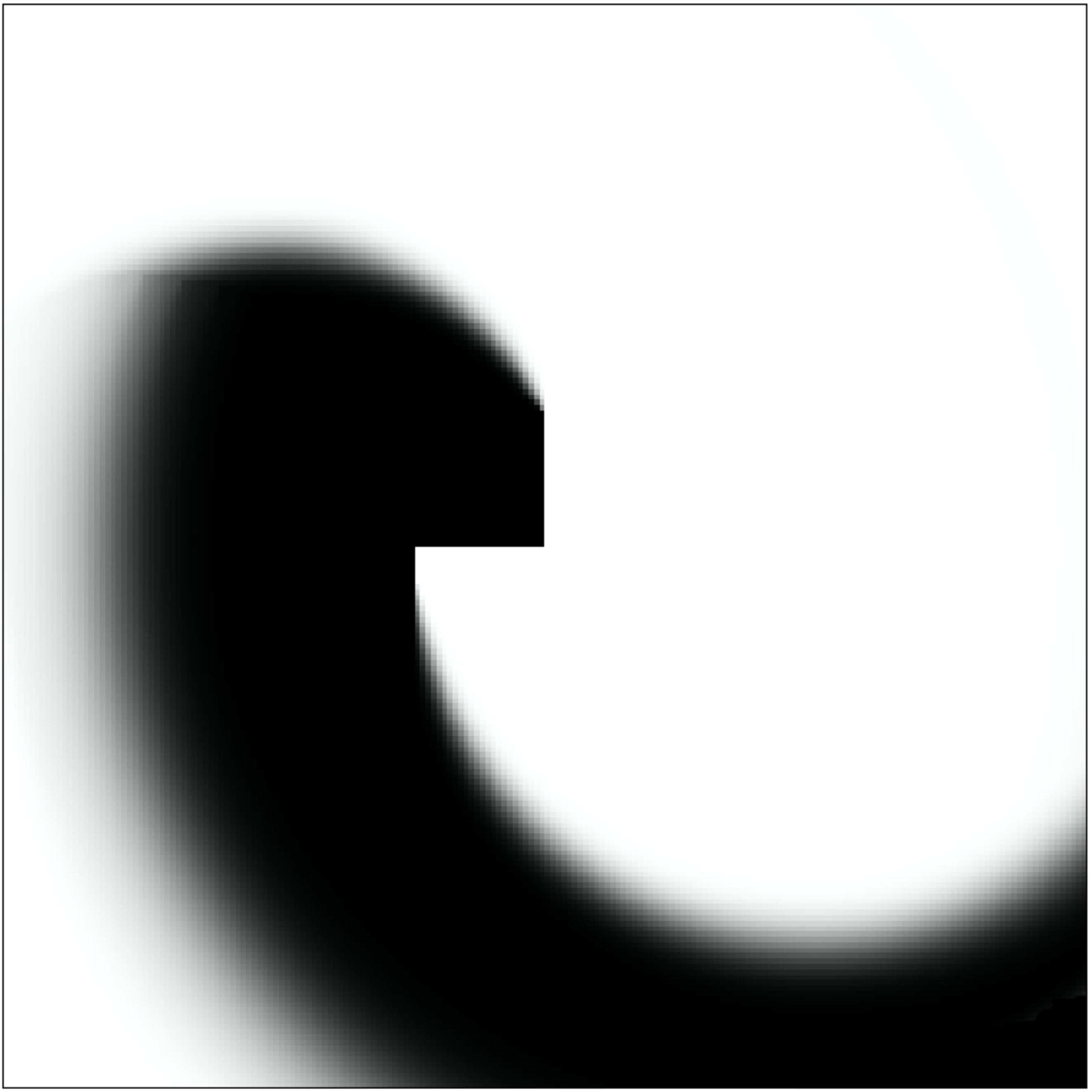} \\
(a) \qquad & \qquad  (b) \\ \\
\includegraphics[width=0.4\textwidth]{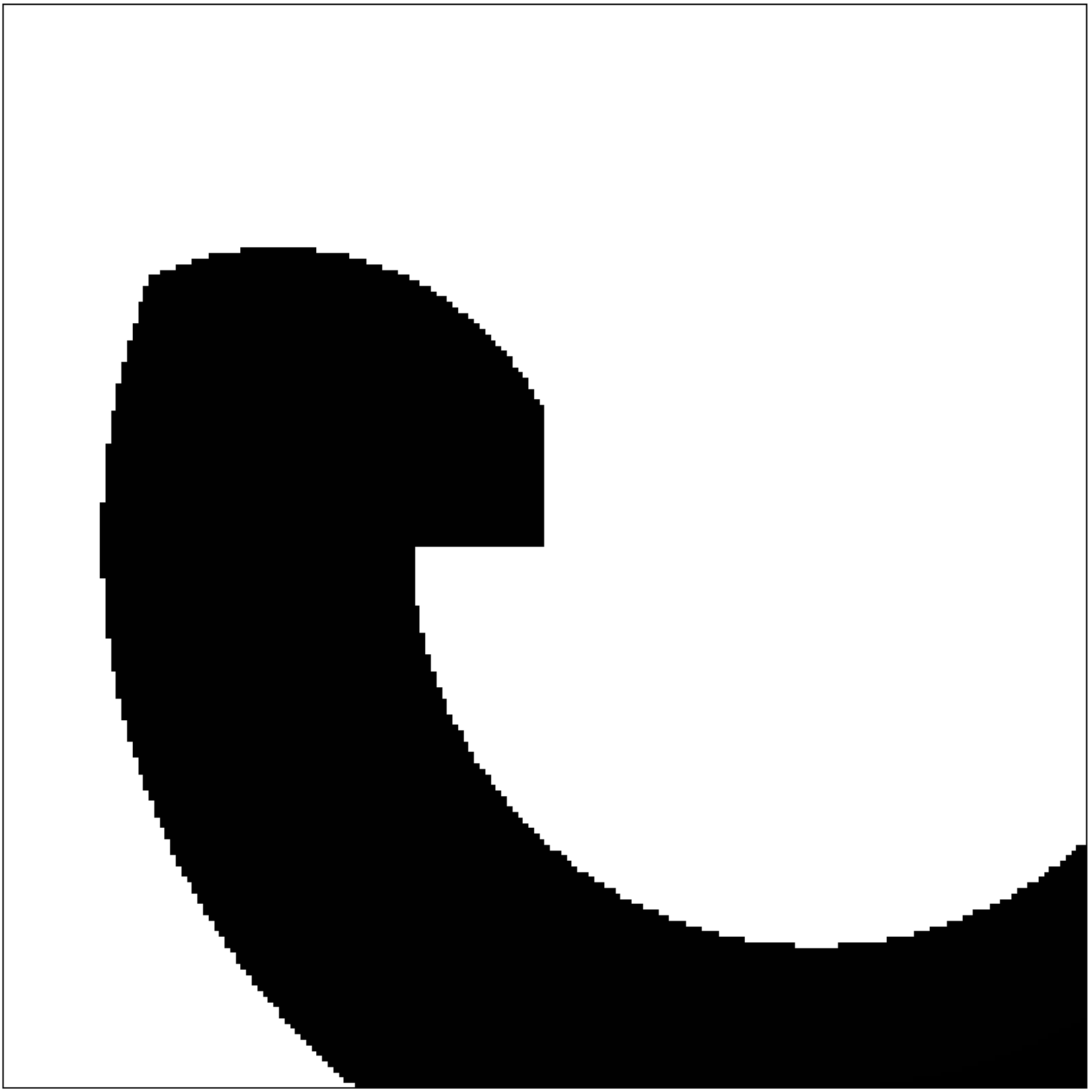} \qquad & \qquad 
\includegraphics[width=0.4\textwidth]{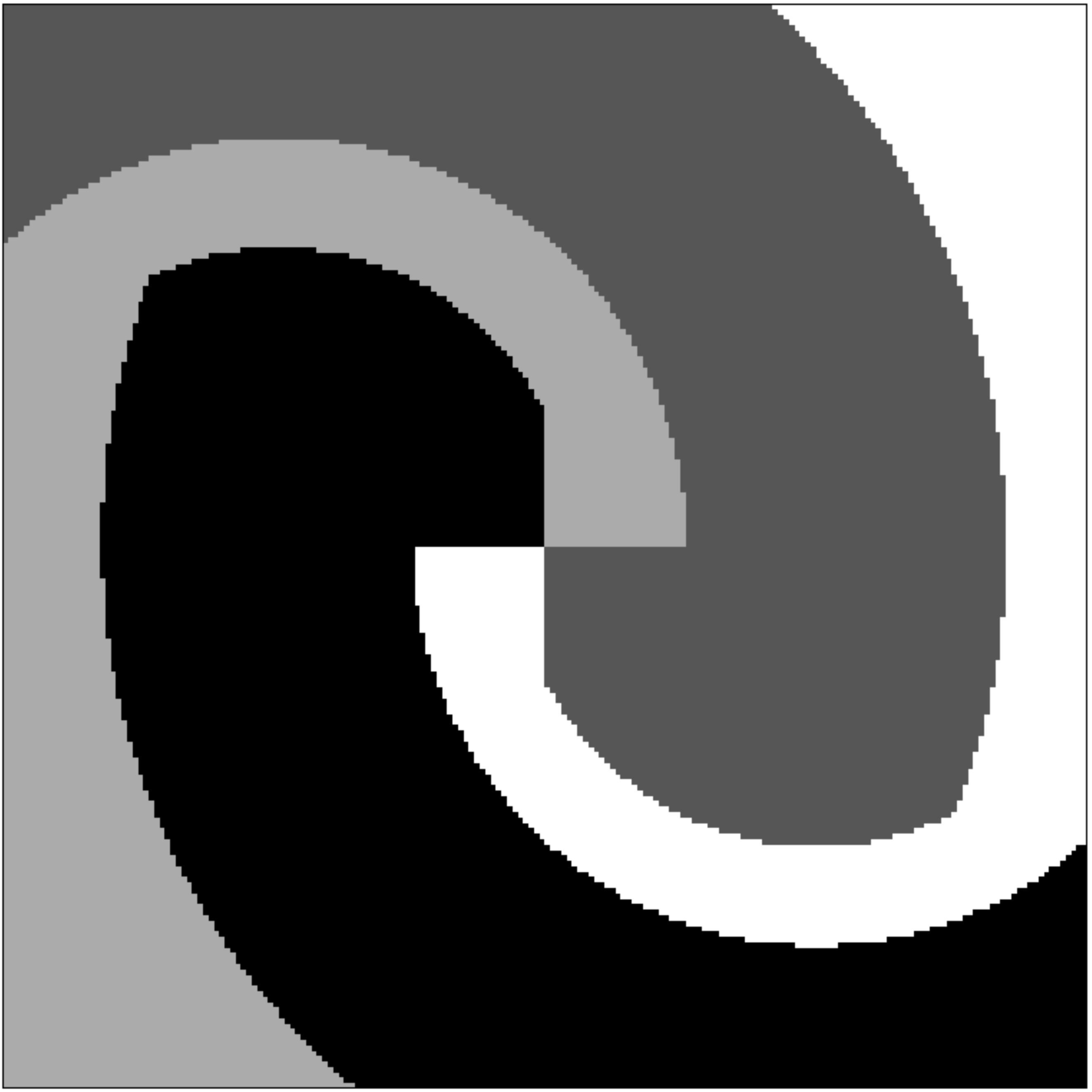} \\
(c) \qquad & \qquad  (d)
\end{tabular}
\caption{Creation of patches for a test dynamics, $R$=4, $\Omega_0$=small ball in the centre: (a) Select a subdomain $\Omega_0^j$ of the target $\Omega_0$. (b) Find the nodes which depend, at least partially, on $\Omega_0^j$. (c) Define $\Omega^j$ projecting the color in a binary value. (d) Assemble all patches.}
\label{fig:creation_patchy}
\end{center}
\end{figure}

\noindent{\bf Patchy Algorithm:}
\begin{enumerate}
\item[Step 1.] (Computation on $\Gcoarse$). We solve the equation on $\Gcoarse$ by means of the classical domain decomposition algorithm described in Section \ref{sec:dd_algorithm}. 
For coherence we choose a (static) decomposition with $R$ subdomains (as the number of patches). This leads to the function $\ucoarse_P$.\\
\item[Step 2.] (Interpolation on $G$). We define the function $U_P^{(0)}$ on the fine grid $G$ by interpolation of the values $\ucoarse_P$. 
Then, we compute the approximate optimal control 
\begin{equation}\label{feedbackcontrol}
a^*_{\widetilde k}(x_i)=\arg\min_{a\in A }\{I[U_P^{(0)}](x_i+h_{i,a}f(x_i,a))+ h_{i,a}\}\,,\qquad x_i\in G.
\end{equation}
Even if $a^*_{\widetilde k}$ is defined on $G$, we still use the subscript $\widetilde k$ to stress that the optimal control is computed using only coarse information. We delete $\Gcoarse$ and $\ucoarse_P$.\\
\item[Step 3.] (Main cycle) For every $j=1,\ldots,R$,\\
\begin{enumerate}
\item[Step 3.1.] (Creation of the $j$-th patch). Using the (coarse) optimal control $a^*_{\widetilde k}$, we find the nodes of the grid $G$ that have the part $\Omega_0^j$ of the target in their numerical domain of dependence. 
This procedure defines the $j$-th patch, naturally following the (approximate) optimal dynamics. This step will be detailed later in this section.\\
\item[Step 3.2.] (Computation in the $j$-th patch). 
We apply iteratively the scheme \eqref{iterativeSLscheme} in the $j$-th patch until convergence. Boundary conditions will be discussed later in this section. 
\end{enumerate}
\medskip
\item[ Step 4.] (Merging) All the solutions are merged together. This leads to the final solution $U_P$.
\end{enumerate}
\medskip
\textit{Details on Step 3.1}. The basic idea we adopt here is dividing the whole domain \textit{starting from a partition of the target only}, and let the dynamics make a partition of the rest of the domain (see \cite{BCGVpp} for a similar idea in the context of parallelization of Fast Marching method for the Eikonal equation).  
More precisely, once the target $\Omega_0$ is divided in $R$ parts, we associate each part to a colour indexed by a number $j=1,\ldots,R$. Assume for instance that $\Omega_0$ is a ball at the center of the domain and focus on the subset of the target with a generic colour $j$, denoted by $\Omega_0^j$, see Fig.\ \ref{fig:creation_patchy}(a). 
The goal is finding the subset of the domain $\Omega$ which has $\Omega_0^j$ as numerical domain of dependence. 
First, we initialize the grid nodes with the values $\phi_i$ as follows:
$$
\phi_i=\left\{
\begin{array}{ll}
1\,, & x_i\in\Omega_0^j\cap G \\
0\,, & x_i\in G\backslash\Omega_0^j
\end{array}
\right.\,,
\qquad
i=1,...,N.
$$
Then, employing the approximation of the optimal control given by $a^*_{\widetilde k}$, we solve the following \emph{ad hoc} discrete equation, 
\begin{equation}\label{HJBridotta_discreta}
\phi_i=I[\phi](x_i+h_if(x_i,a^*_{\widetilde k}(x_i)))\,,\qquad i=1,...,N,
\end{equation}
which is similar to the fixed-point scheme \eqref{SLscheme} for the main equation.
Here $h_i>0$ is chosen in such a way that $|h_i f(x_i,a^*_{\widetilde k}(x_i))| = k$. 
Once the computation is completed, the whole domain will be divided in three zones:
$$
\Lambda^j_1=\{x_i:\phi_i=1\}\,,\quad \Lambda^j_2=\{x_i:\phi_i=0\}\,, \quad \Lambda^j_3=\{x_i:\phi_i\in(0,1)\}\,,
$$ 
see Fig.\ \ref{fig:creation_patchy}(b). 
Note that $\Lambda^j_3$ will be nonempty because the interpolation operator $I$ in the scheme \eqref{HJBridotta_discreta} mixes the values $\phi_i$ through a convex 
combination, thus producing values in $[0,1]$ even if the initial datum is in $\{0,1\}$. Since we need a sharp division of the domain, we ``project'' the colour $j$ into a binary value
\begin{equation}\label{proiezionecolore}
\widehat\phi_i=\left\{
\begin{array}{ll}
1\,, & \phi_i\geq \frac12 \\
0\,, & \phi_i<\frac12 
\end{array}
\right.\,,
\qquad
i=1,...,N
\end{equation}
and then we define the subdomain $\Omega^j=\{x_i\in G\backslash\Omega_0^j : \widehat\phi_i=1\}$ as the $j$-th patch, see Fig.\ \ref{fig:creation_patchy}(c). 
Once all the patches $j=1,\ldots,R$ are computed, they are assembled together on the grid $G$. Thus the grid results to be divided into $R$ patches, each associated to a different colour, as shown in Fig.\ \ref{fig:creation_patchy}(d). 

The main point here is that patches $\Omega^j$'s are constructed to be invariant with respect to the optimal dynamics, meaning that the solution of the equation in each patch will not depend on the solution in other patches. 
This is equivalent to state that there is no crossing information through the boundaries of the patches. 

We stress that Step 3.1 of the algorithm is not expensive, even if it is performed on the fine grid $G$. 
The reason for that is the employment of the pre-computed optimal control $a^*_{\widetilde k}$ in the equation (\ref{HJBridotta_discreta}), which avoids the evaluation of the minimum (see the scheme (\ref{iterativeSLscheme})). Moreover,
the stop criterion for the fixed-point iterations used to solve (\ref{HJBridotta_discreta}) can be very rough, since we project the colors at the end and then we do not need precise values. 

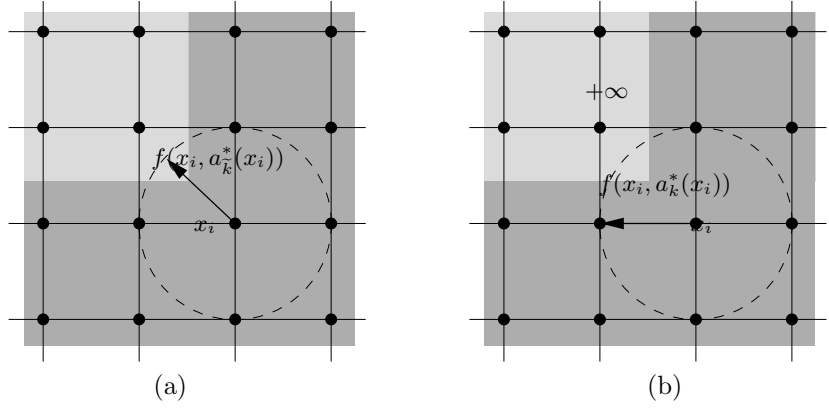
\begin{figure}[h!]
\begin{center}
\begin{tabular}{c}
\includegraphics[width=.9\textwidth]{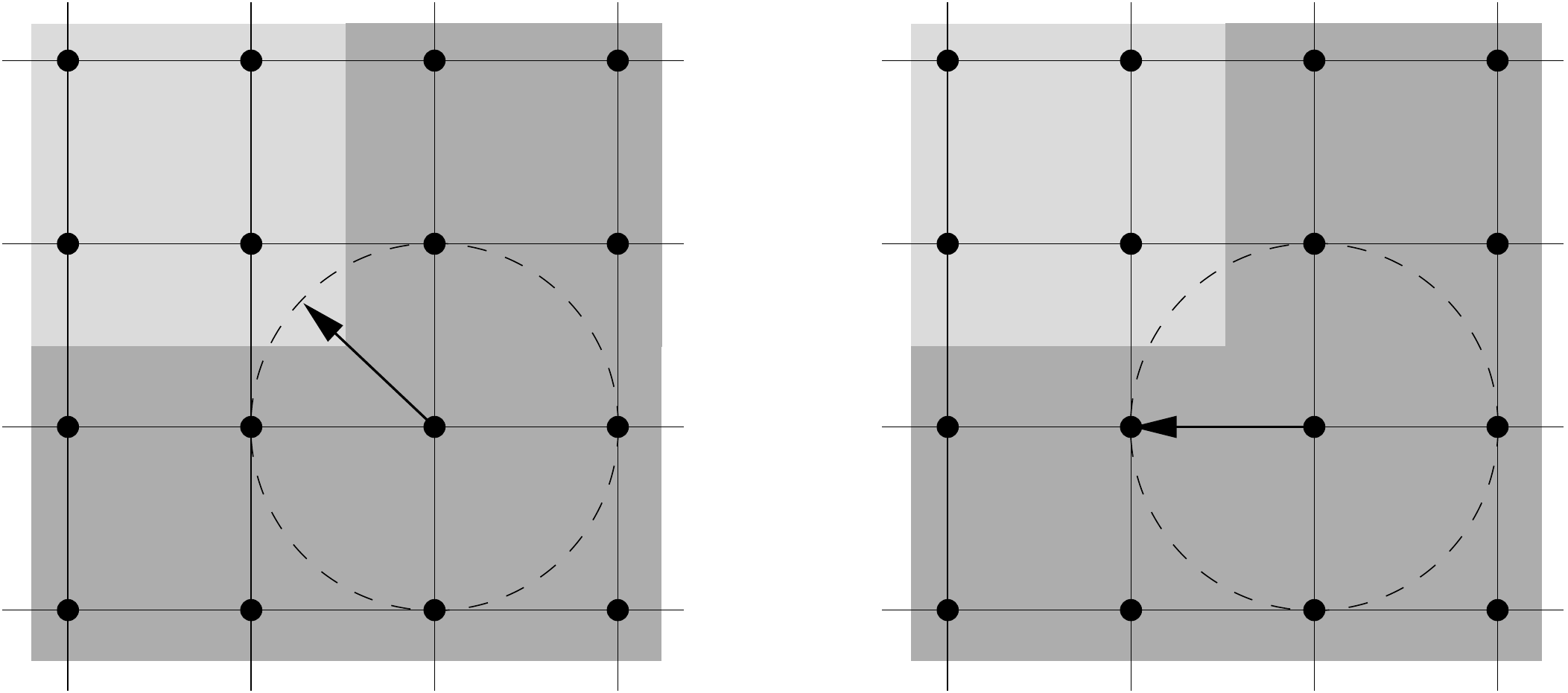}  
\put(-256,75){\small $f(x_i,a^*_{\widetilde k}(x_i))$}
\put(-240,50){\small $x_i$}
\put(-92,100){\small $+\infty$}
\put(-87,65){\small $f(x_i,a^*_k(x_i))$}
\put(-52,50){\small $x_i$}
\\
(a) \hskip6cm  (b) 
\end{tabular}
\caption{(a) Discretization vs decomposition invariance: the node $x_i$ belongs to the dark-gray patch, but it is influenced by the light-gray patch (through the coarse control $a^*_{\widetilde k}(x_i)$) . (b) 
The cure: state constraint boundary conditions force the patches to be completely independent, but change the optimal vector field $f(x_i,a^*_k(x_i))$ producing an error in the final solution.}
\label{fig:independence}
\end{center}
\end{figure}

\medskip
\textit{Details on boundary conditions}. In Step 3.2 of the algorithm,
the computation of the value function is performed independently in
each patch, thus we have to impose boundary conditions on the
boundaries of the patches. A natural choice is the employment of
$U_P^{(0)}$ (obtained in Step 2) as Dirichlet boundary condition,
but in some early tests \cite{P11} we observed that this choice leads to
reasonable results for $U_P$ even if the domain decomposition
is completely incorrect (i.e.\ patches are not at all invariant).
More precisely, if the decomposition is not invariant, the accuracy of
the final solution $U_P$ is comparable to that of $U_P^{(0)}$,
otherwise it is comparable to that of $U_{\!D\!D}$ computed on the
same fine grid (which is the best one can do). This point will be
discussed later in Section \ref{sec:patchyerror}.
As a first study we prefer to impose a boundary condition which does
not require any a-priori information on the solution outside the
patches,
in order to check if they are genuinely independent. This motivated us
to use state constraint boundary conditions, which force the optimal
direction $f(x_i,a^*_k)$ to point inside
the patch (see Fig.\ \ref{fig:independence} and its caption). This
choice produces an error that can be evaluated comparing $U_P$ with
$U_{\!D\!D}$ and that we consider as a degree of invariance of the patchy decomposition.

Let us remark that once we have shown that patches are actually (almost) invariant, we \emph{can} impose Dirichlet boundary conditions, further improving the quality of the solution.
This also allows one to handle the case of dynamics such that state constraints condition cannot be satisfied everywhere (i.e.\ there is no control allowing the state to remain in the patch), as in the example discussed in Section \ref{krener}.

\medskip

\textit{Remark 3.1 (patches as a partition of $G$).} We have no guarantee that patches $\Omega^j$'s do not overlap or that they cover the whole domain. On the overlapping zones we can simply choose a colour at random. Instead, if they do not cover all the domain we can repeat the computation in the not-coloured nodes relaxing the condition in \eqref{proiezionecolore}, i.e.\ choosing a different value for $1/2$. Alternatively, in the case of isolated not-coloured points, we can assign to them the colour of their neighbours.


\medskip

\subsubsection*{How to parallelize the algorithm} 
The patchy algorithm can be parallelized in two ways.

\begin{itemize}
 \item \emph{Method 1}: Patches are processed \textit{one after the other} and the computation in each patch is parallel, assigning a batch of nodes to each processor. \\
\item \emph{Method 2}: Patches are \textit{distributed among processors} and the computation of each patch is serial. 
\end{itemize}
\medskip
The first strategy is designed for shared-memory architectures and gives priority to \textit{saving CPU time}, while the second strategy is designed for distributed-memory architectures and gives priority to \textit{saving memory allocation}. The difference in CPU time comes from the fact that, using the first method, processors are active all the time, while, using the second method, it can happen that one processor finishes its jobs and there are no more patches to be computed, so it remains idle. 
We stress again that, employing Method 2, the independence of the patches allows processors to not communicate until the end of their task, saving heavy overhead in distributed-memory architectures. 

All tests presented in this paper are performed implementing Method 1 on a shared-memory architecture, and in the following we will always refer to this choice.
\section{Numerical investigation in dimension two}
In this section we first list the dynamics considered for the numerical tests. Then, we investigate the optimality of the patchy decomposition and the performance of the algorithm with respect to the classical domain decomposition.

Numerical tests were performed on a server Supermicro 8045C-3RB using 1 CPU Intel Xeon Quad-Core E7330 2.4 Ghz with 32 GB RAM, running under Linux Gentoo operative system.
\subsection{Choice of benchmarks}\label{benchmarks2d}
\begin{figure}[h!]
\begin{center}
\begin{tabular}{cc}
\includegraphics[width=0.49\textwidth]{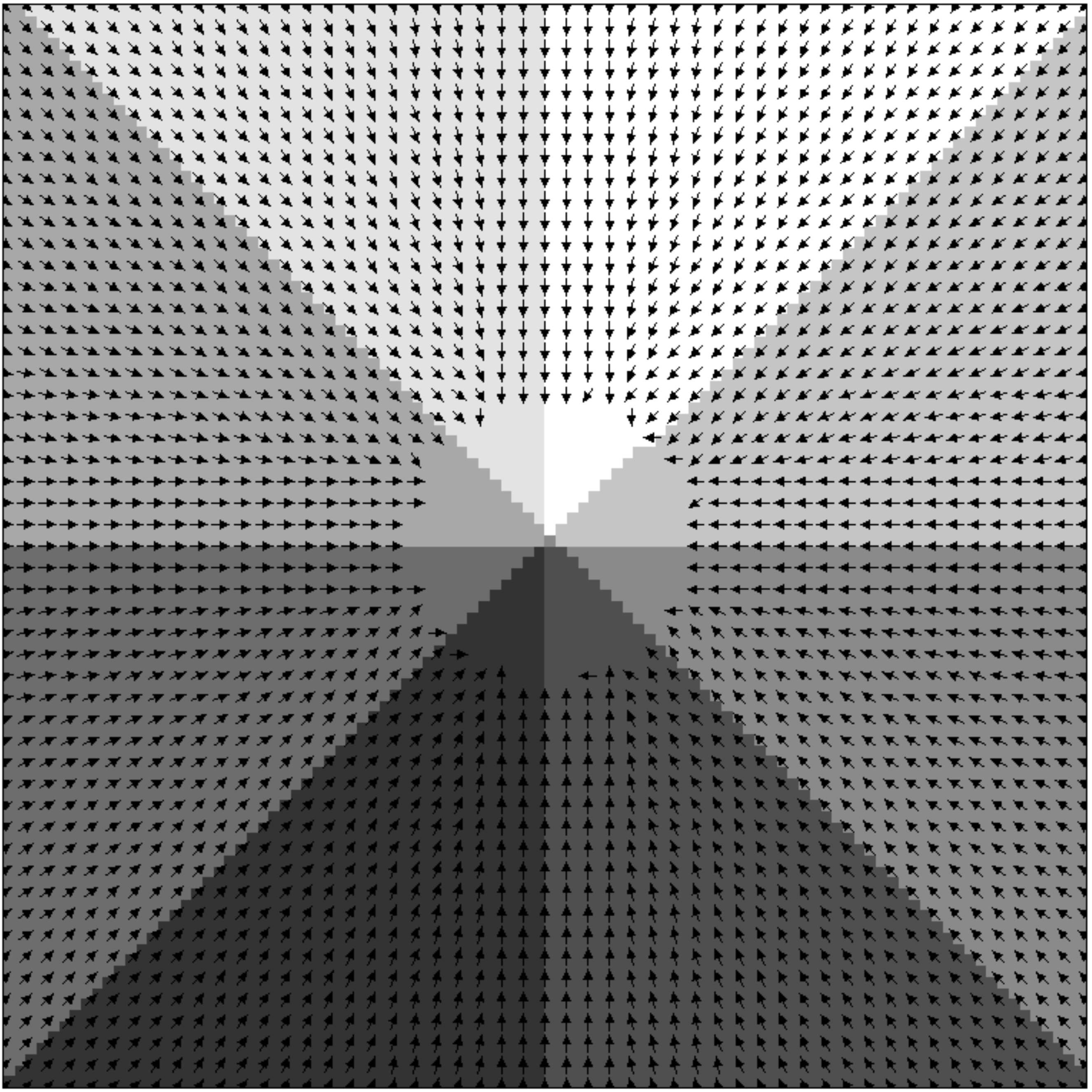} &
\includegraphics[width=0.49\textwidth]{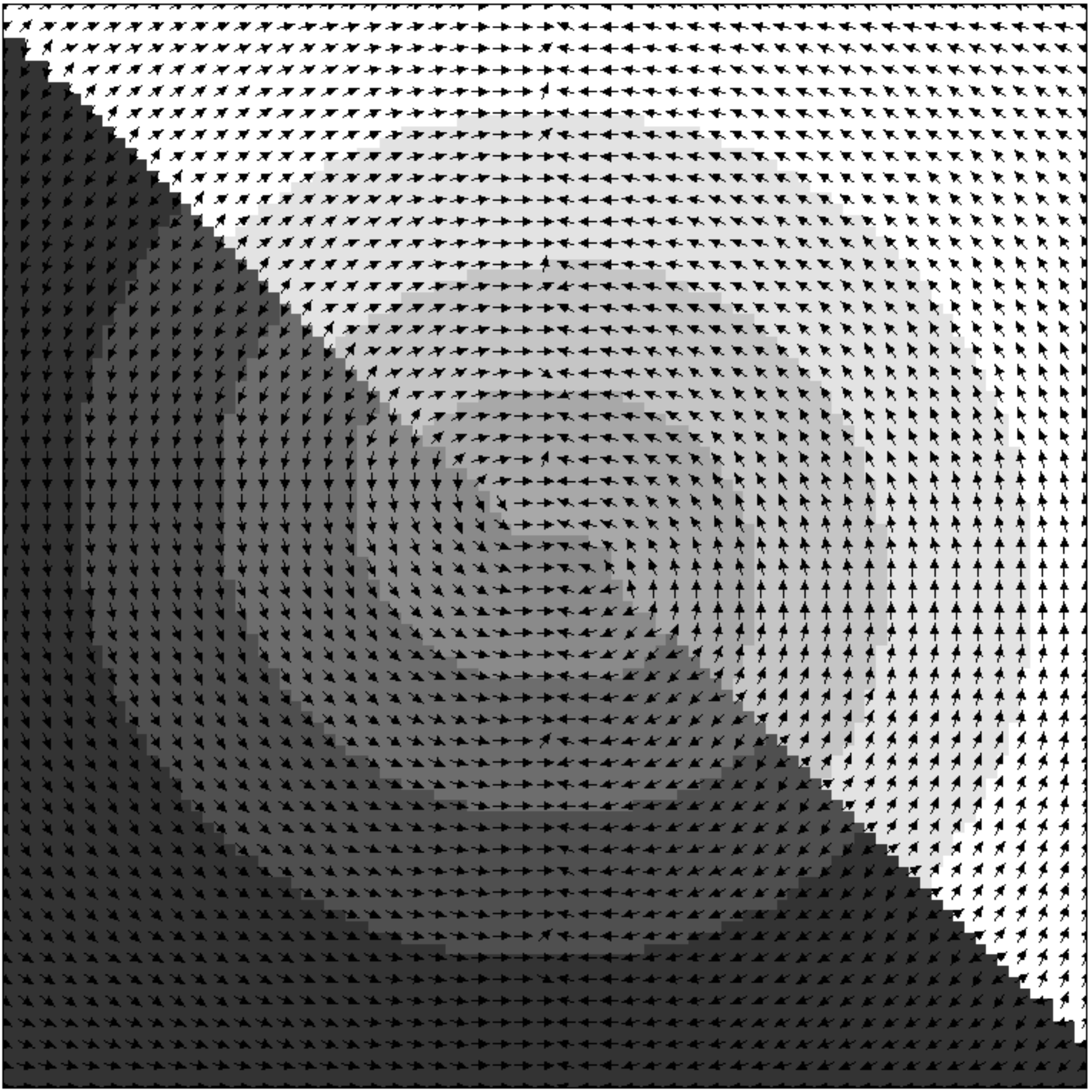} \\
(a)  &  (b)
\end{tabular}
\begin{tabular}{cc}
\includegraphics[width=0.49\textwidth]{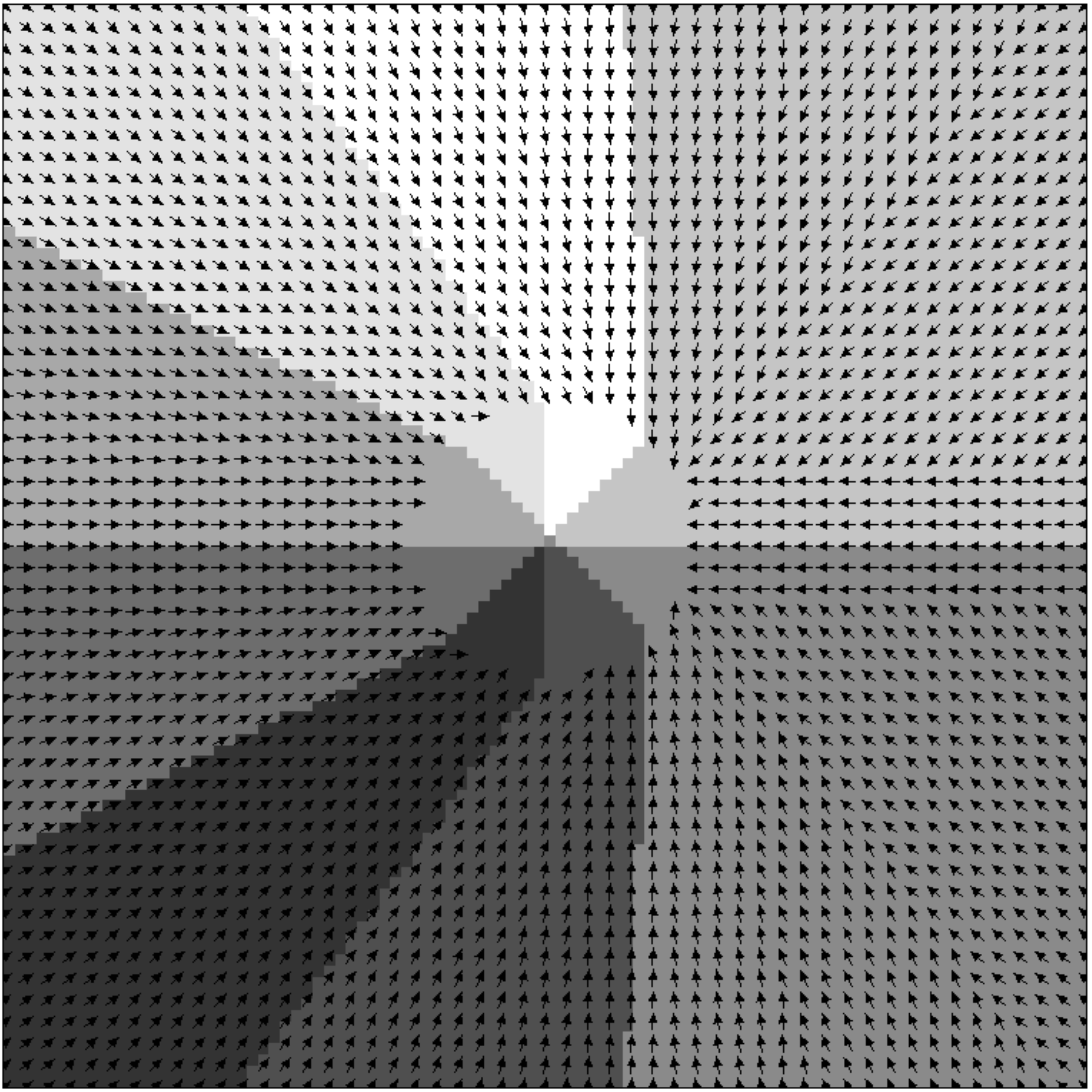} &
\includegraphics[width=0.49\textwidth]{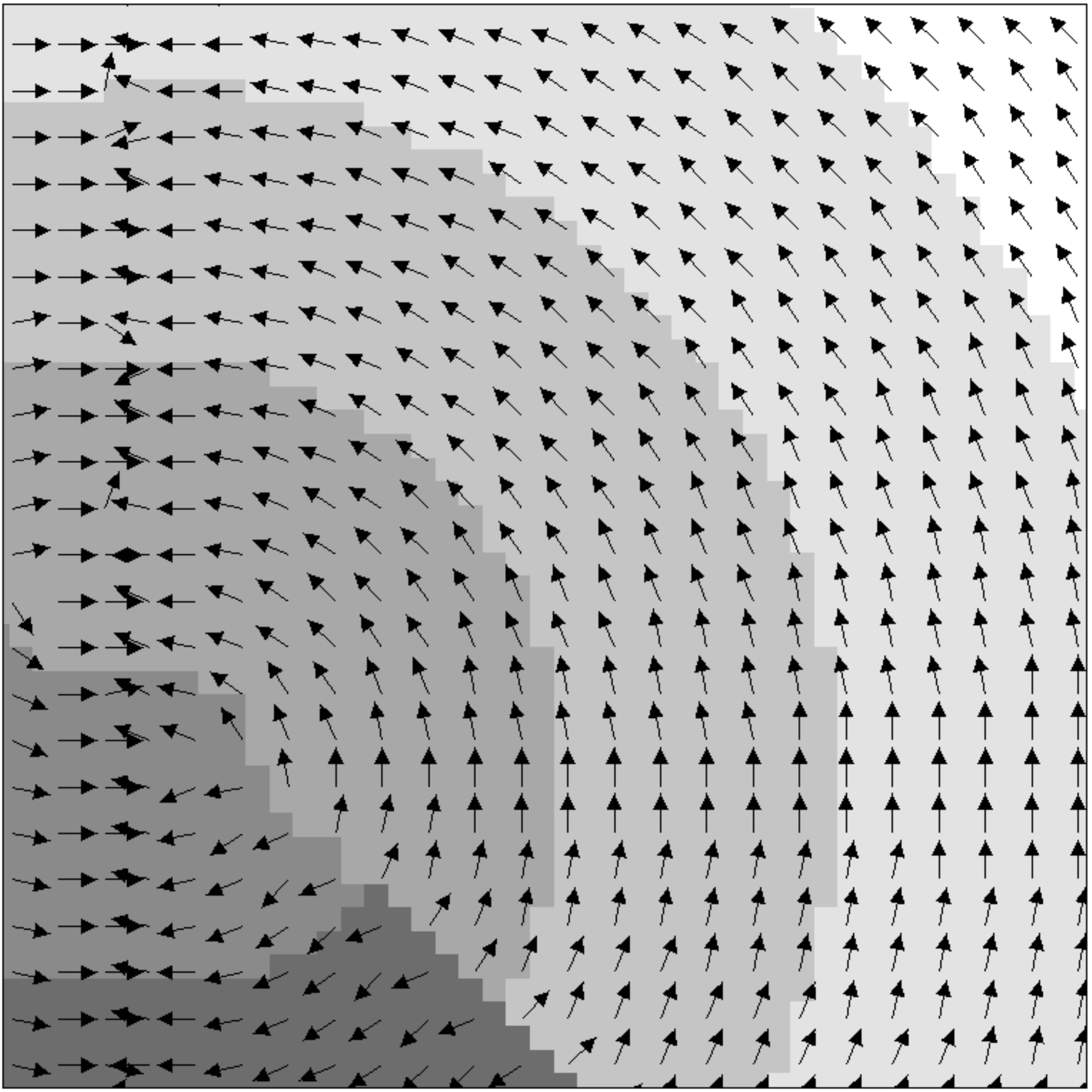} \\
(c)  &  (d)
\end{tabular}
\caption{Patchy decompositions with $R=8$, $N_c=32$, $\Ncoarse=50^2$ and $N=100^2$. For visualization purposes not all the arrows are shown. (a) Eikonal, (b) Fan, (c) Zermelo, (d) a detail of Fan.}
\label{fig:patchy_decomposition}
\end{center}
\end{figure}

We will test the method described above against three minimum time problems of the form (\ref{HJB}). They are listed in Table \ref{tab:test2d}. The numerical domain is $\Omega=[-2,2]^2$ for all tests.
\begin{table}[h!]
\caption{Two-dimensional numerical tests}
\label{tab:test2d}
\begin{center}
\begin{tabular}{|l|l|l|l|l|}
\hline Name    & $d$ & $f(x_1,x_2,a)$   & $A$ & $\Omega_0$ \\ \hline
\hline Eikonal & 2   & $a$              & $B_2(0,1)$ & $B_2(0,0.5)$ \\
\hline Fan     & 2   & $|x_1+x_2+0.1|a$ & $B_2(0,1)$ & $\{x_1=0\}$ \\
\hline Zermelo & 2   & $2.1a+(2,0)$     & $B_2(0,1)$ & $B_2(0,0.5)$ \\
\hline
\end{tabular}
\end{center}
\end{table}
In Fig.\ \ref{fig:patchy_decomposition} we show 
the patchy decomposition for the three dynamics described above in the case $R=8$, $\Ncoarse=50^2$ and $N=100^2$. The number of points used for discretizing $A$ is $N_c=32$. We also superimpose the coarse optimal vector field $f(x,a^*_{\widetilde k})$ to show that patches are (almost) invariant with respect to the optimal dynamics. Indeed, only a few arrows cross from a patch to another, which is exactly the case discussed in Fig.\ \ref{fig:independence}(a).  
Note that patches cover the whole domain but in general they are not equivalent in terms of area, even if the target $\Omega_0$ was divided in $R=8$ equal parts to generate the decomposition.
\subsection{Optimality of the patchy decomposition}\label{sec:patchyerror}
In this section we compare the solution $U_P$ of the patchy algorithm with that of the classical domain decomposition method $U_{D\! D}$, both computed on the same fine grid 
by means of the 
scheme \eqref{iterativeSLscheme}. Let us denote by $E$ the difference 
$$
E:=U_{\!P}-U_{\!D\!D}
$$ 
that in the following will be referred to as \textit{patchy error}. 
In particular we study the quantities 
$$E_1:=\frac1N\sum_{i=1}^N |E_i|
\qquad\mbox{and}\qquad 
E_\infty:=\max_{i=1,...,N}|E_i|$$ 
as the grid nodes $\Ncoarse$ and $N$ change. Error $E$ is exclusively due to the fact that patches are not completely dynamics-invariant and then it will be considered as a degree of the invariance of the patchy decomposition. Let us stress that we employ state constraint boundary conditions, 
as discussed in the previous section (see Fig.\ \ref{fig:independence}(b)).

We report the results for $R=16$, which is the largest number of patches and also the worst case we tested. Indeed, error $E$ necessarily increases as $R$ increases because the number of boundaries increases. 
Results for the three dynamics are shown in Tables \ref{tab:patchy_error_eikonal}, \ref{tab:patchy_error_fan} and \ref{tab:patchy_error_zermelo}. 

\begin{table}[h!]
\small
\caption{Patchy error $E_1$ ($E_\infty$). Dynamics: Eikonal, $N_c=32$, $R=16$}
\label{tab:patchy_error_eikonal}
\begin{center}
\begin{tabular}{|l|c|c|c|c|c|}
\hline $ $ & $N=50^2$ & $N=100^2$ & $N=200^2$ & $N=400^2$ & $N=800^2$ \\
\hline
\hline \Ncoarsebis =$50^2$   & $\ba 0.02725\\(0.960)\ea$ & $\ba 0.01719\\(1.856)\ea$  & $\ba 0.00637\\(0.048)\ea$  & $\ba 0.00406\\(0.034)\ea$  & $\ba 0.00300\\(0.026)\ea$ \\
\hline \Ncoarsebis =$100^2$  &      --                   & $\ba 0.00550\\(0.046)\ea$  & $\ba 0.00181\\(0.023)\ea$  & $\ba 0.00087\\(0.042)\ea$  & $\ba 0.00031\\(0.008)\ea$ \\
\hline \Ncoarsebis =$200^2$  &      --                   &      --                    & $\ba 0.00237\\(0.029)\ea$  & $\ba 0.00075\\(0.013)\ea$  & $\ba 0.00025\\(0.008)\ea$ \\
\hline \Ncoarsebis =$400^2$  &      --                   &      --                    &      --                    & $\ba 0.00069\\(0.016)\ea$  & $\ba 0.00037\\(0.010)\ea$ \\
\hline \Ncoarsebis =$800^2$  &      --                   &      --                    &      --                    &      --                    & $\ba 0.00025\\(0.008)\ea$ \\
\hline
\end{tabular}
\end{center}
\end{table}

\begin{table}[h!]
\small
\caption{Patchy error  $E_1$ ($E_\infty$). Dynamics: Fan, $N_c=32$, $R=16$}
\label{tab:patchy_error_fan}
\begin{center}
\begin{tabular}{|l|c|c|c|c|c|}
\hline $ $ & $N=50^2$ & $N=100^2$ & $N=200^2$ & $N=400^2$ & $N=800^2$ \\
\hline
\hline \Ncoarsebis =$50^2$   & $\ba 0.08706\\(3.023)\ea$ & $\ba 0.00769\\(1.507)\ea$  & $\ba 0.00231\\(0.315)\ea$  & $\ba 0.00106\\(0.263)\ea$  & $\ba 0.00069\\(0.263)\ea$ \\
\hline \Ncoarsebis =$100^2$  &      --                   & $\ba 0.00712\\(1.502)\ea$  & $\ba 0.00200\\(0.149)\ea$  & $\ba 0.00069\\(0.095)\ea$  & $\ba 0.00037\\(0.095)\ea$ \\
\hline \Ncoarsebis =$200^2$  &      --                   &      --                    & $\ba 0.00200\\(0.111)\ea$  & $\ba 0.00069\\(0.061)\ea$  & $\ba 0.00025\\(0.037)\ea$ \\
\hline \Ncoarsebis =$400^2$  &      --                   &      --                    &      --                    & $\ba 0.00069\\(0.079)\ea$  & $\ba 0.00025\\(0.037)\ea$ \\
\hline \Ncoarsebis =$800^2$  &      --                   &      --                    &      --                    &      --                    & $\ba 0.00025\\(0.037)\ea$ \\
\hline
\end{tabular}
\end{center}
\end{table}

\begin{table}[h!]
\small
\caption{Patchy error  $E_1$ ($E_\infty$). Dynamics: Zermelo, $N_c=32$, $R=16$}
\label{tab:patchy_error_zermelo}
\begin{center}
\begin{tabular}{|l|c|c|c|c|c|}
\hline $ $ & $N=50^2$ & $N=100^2$ & $N=200^2$ & $N=400^2$ & $N=800^2$ \\
\hline
\hline \Ncoarsebis =$50^2$   & $\ba 0.01069\\(0.293)\ea$ & $\ba 0.00994\\(0.059)\ea$  & $\ba 0.00606\\(0.057)\ea$  & $\ba 0.00162\\(0.027)\ea$  & $\ba 0.00037\\(0.016)\ea$ \\
\hline \Ncoarsebis =$100^2$  &      --                   & $\ba 0.00631\\(0.063)\ea$  & $\ba 0.00206\\(0.041)\ea$  & $\ba 0.00069\\(0.023)\ea$  & $\ba 0.00025\\(0.016)\ea$ \\
\hline \Ncoarsebis =$200^2$  &      --                   &      --                    & $\ba 0.00244\\(0.039)\ea$  & $\ba 0.00075\\(0.023)\ea$  & $\ba 0.00025\\(0.016)\ea$ \\
\hline \Ncoarsebis =$400^2$  &      --                   &      --                    &      --                    & $\ba 0.00069\\(0.020)\ea$  & $\ba 0.00031\\(0.015)\ea$ \\
\hline \Ncoarsebis =$800^2$  &      --                   &      --                    &      --                    &      --                    & $\ba 0.00025\\(0.016)\ea$ \\
\hline
\end{tabular}
\end{center}
\end{table}

We see that the first line of each table ($\Ncoarse$=$50^2$) reports in some cases unsatisfactory results, caused by the excessive roughness of the grid $\Gcoarse$. 
Even the case $\Ncoarse$=$N$=$50^2$ (i.e.\ the grid is not refined at all) is not satisfactory. This can be explained by recalling that, even if $\Ncoarse$=$N$, the computations on the two grids are not identical because the second one employs state constraint boundary conditions. 
In the other cases, the behaviour of the error is very good because \textit{it decreases as $N$ increases} (for any fixed $\Ncoarse$). As pointed out in Section \ref{sec:patchy}, if the patches were not invariant with respect to the dynamics, as in the classical domain decomposition algorithm, 
we would not expect such a behaviour here, because of the missing information across the patches. Thus, this shows that our patches are actually (almost) independent.

Tables also show that the $E_\infty$ is always larger than $E_1$, meaning that the error is concentrated in small regions. Indeed, quite often we find a very small number of nodes with a large error near the boundaries of the patches, especially at those nodes where two patches and the target meet. This mainly affects $E_\infty$ but not $E_1$. 
Finally we note that the results are similar for the three dynamics, showing a good robustness even for highly rotating vector fields like that of Fan dynamics.

\begin{figure}[h!]
\begin{center}
\begin{tabular}{cc}
\includegraphics[width=0.53\textwidth]{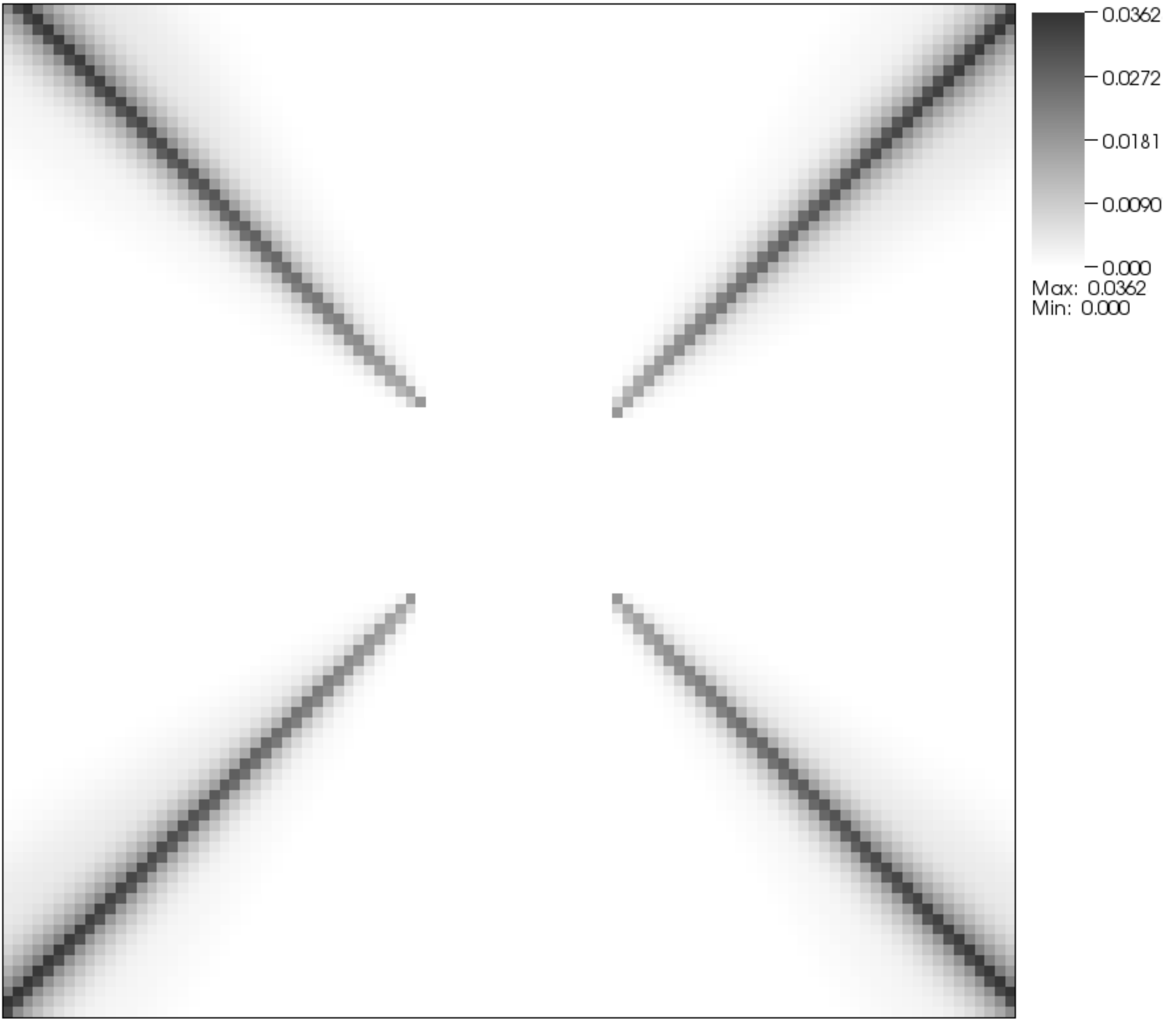} &
\includegraphics[width=0.53\textwidth]{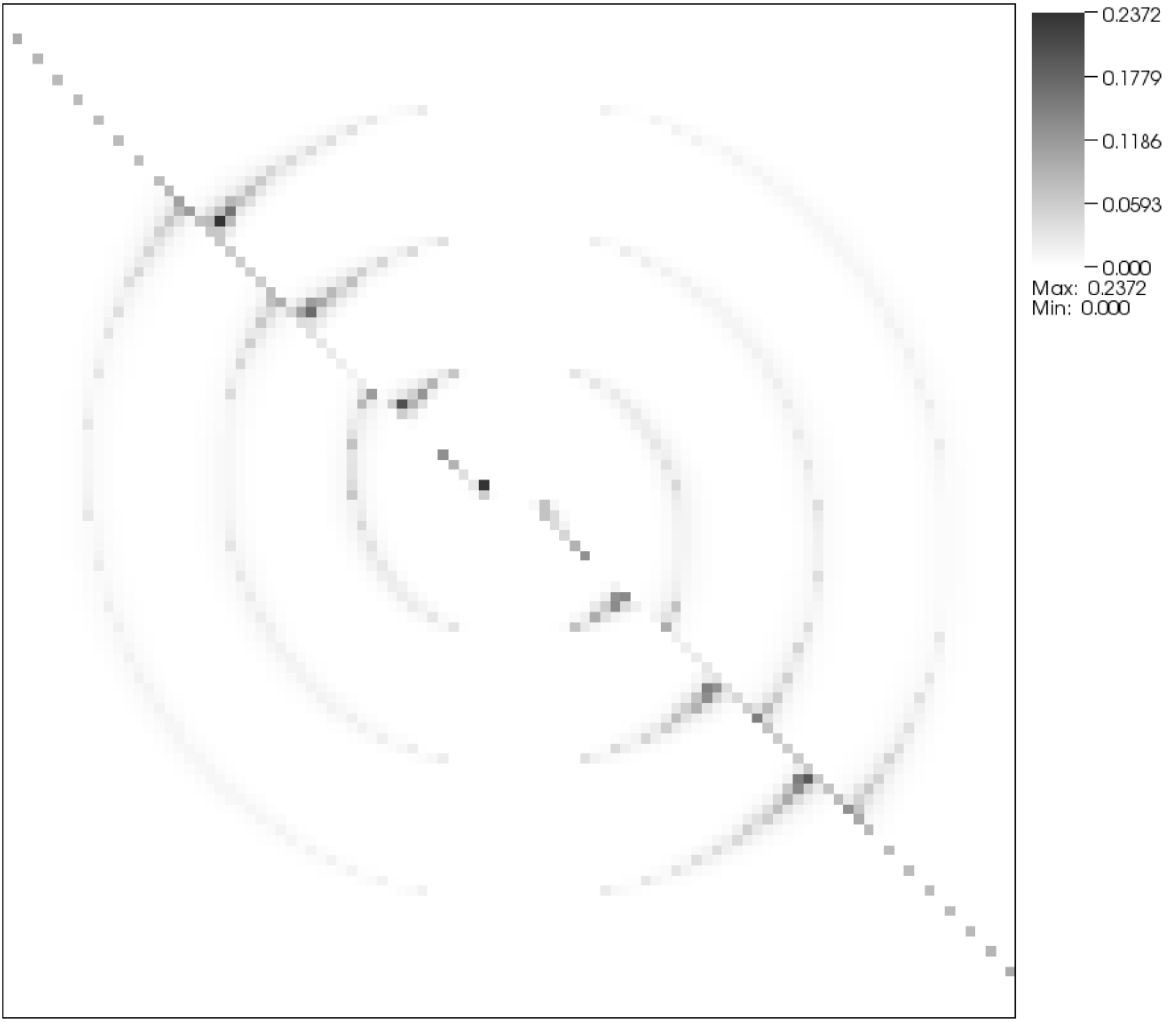} \\ 
(a)  & \hskip-25pt (b) \\
\includegraphics[width=0.53\textwidth]{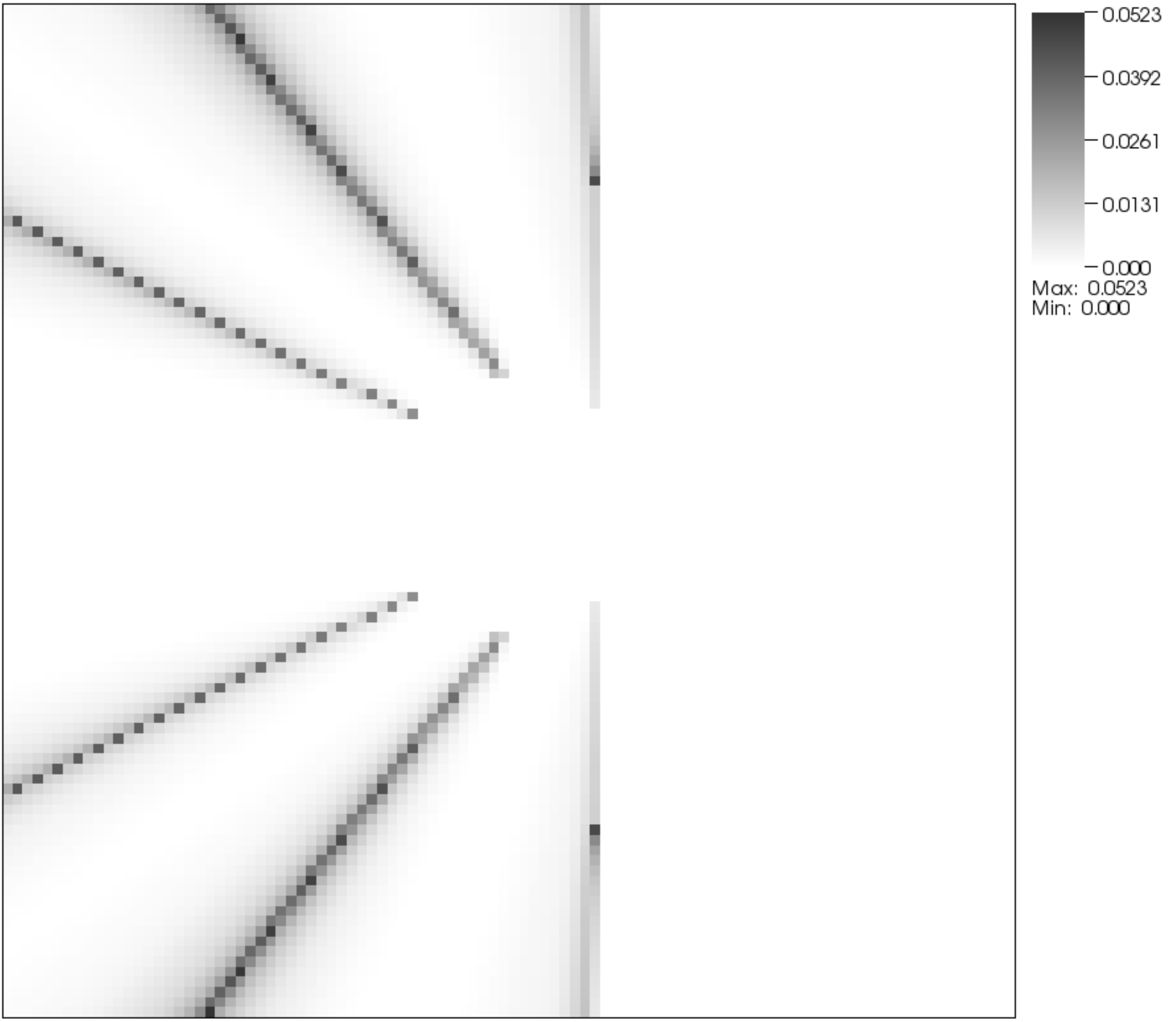} &
\hskip-25pt \includegraphics[width=0.47\textwidth]{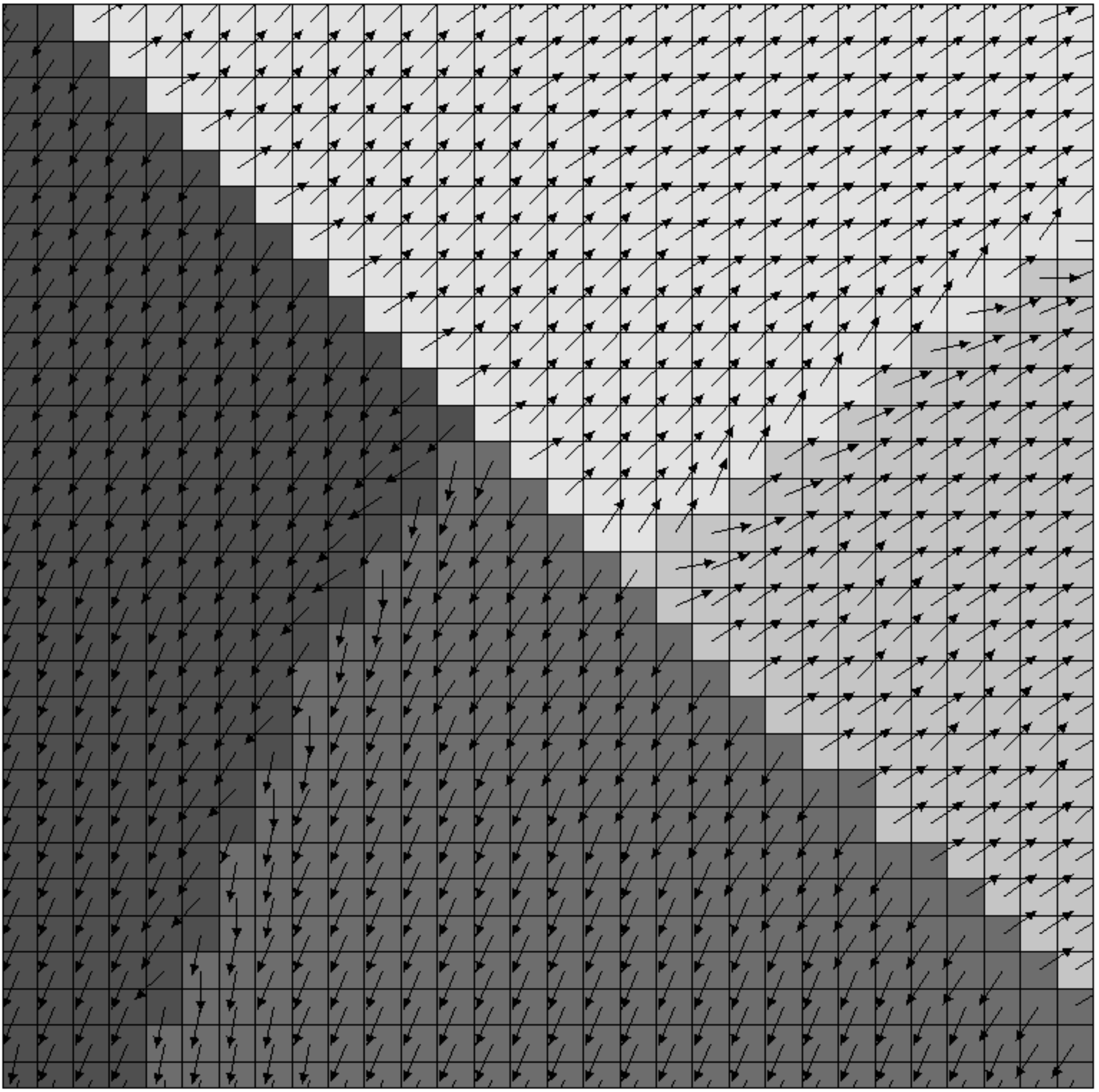} \\
(c) & \hskip-25pt (d) \\
\end{tabular}
\caption{Patchy error $E$, $\Ncoarse=50^2\rightarrow N=100^2$ for (a) Eikonal, (b) Fan, (c) Zermelo. In (d) it is shown a detail of the patchy decomposition for the Fan dynamics, together with the optimal vector field $f(x,a^*_k)$ computed by means of the final patchy solution $U_P$.}
\label{fig:patchy_errorLinf}
\end{center}
\end{figure}

Fig.\ \ref{fig:patchy_errorLinf}-(a,b,c) reports the function $E$ for the three tests, showing one of the most interesting features of the new method, i.e.\ the patchy error 
is concentrated along the boundaries of the patches and does not propagate in the interior. 
In the Eikonal and Zermelo case the error starts from the target and increases as long as characteristics go away. 
In the Fan case, instead, the largest error is found where patches and target meet. Note that in the Eikonal case (a) no error is found where patches boundaries are aligned to the grid, since the optimal direction naturally points inside the patch and the state constraint boundary condition has no effect.
Fig.\ \ref{fig:patchy_errorLinf}-(d) shows a detail of the Fan decomposition along with the approximate optimal vector field computed by means of the final solution $U_P$. The effect of the state constraints is perfectly visible (arrows point unnaturally inward) confirming the fact that each patch is computed independently.

\begin{figure}[h!]
\begin{center}
\begin{tabular}{cc}
\includegraphics[width=0.49\textwidth]{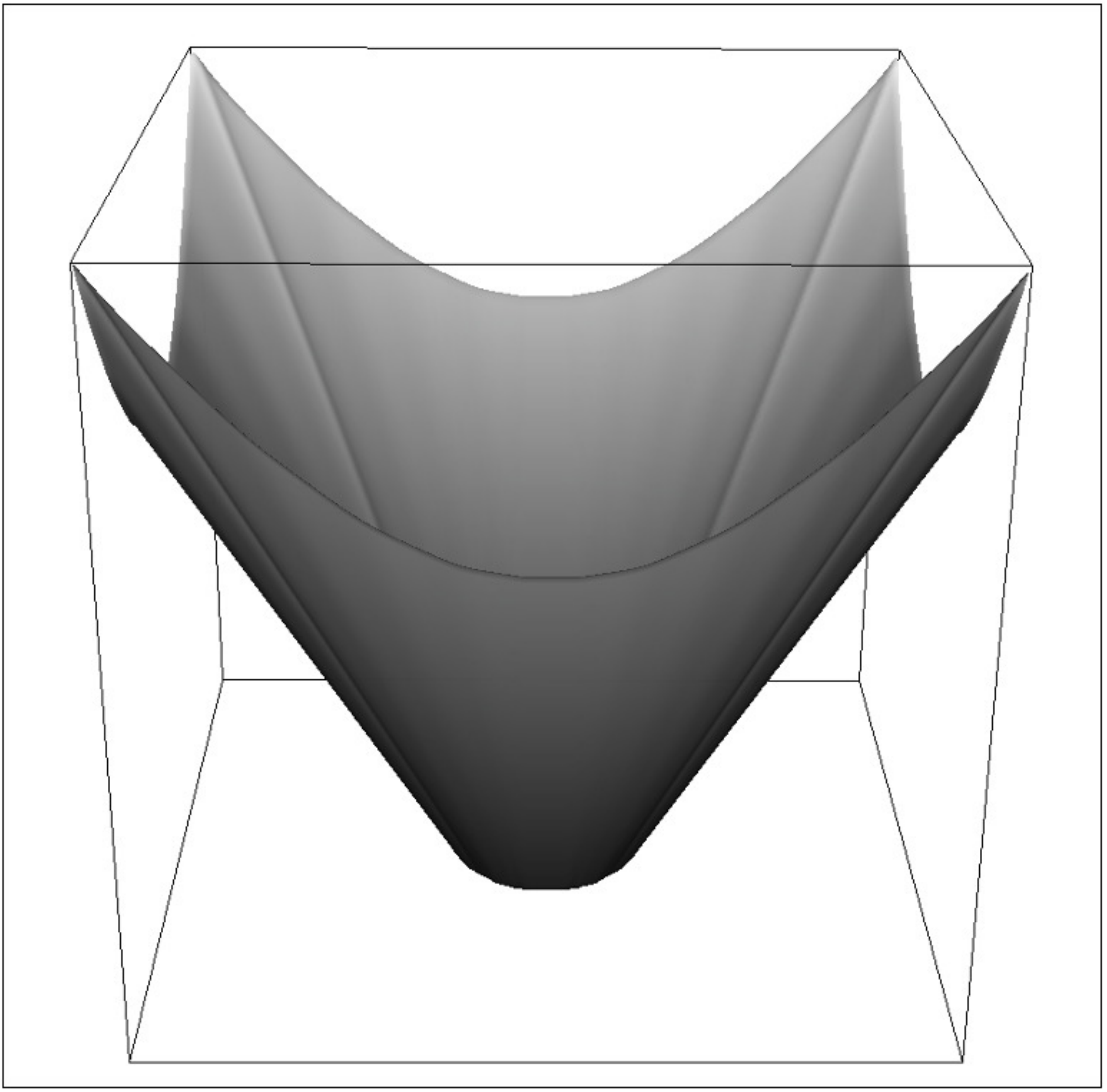} &
\includegraphics[width=0.49\textwidth]{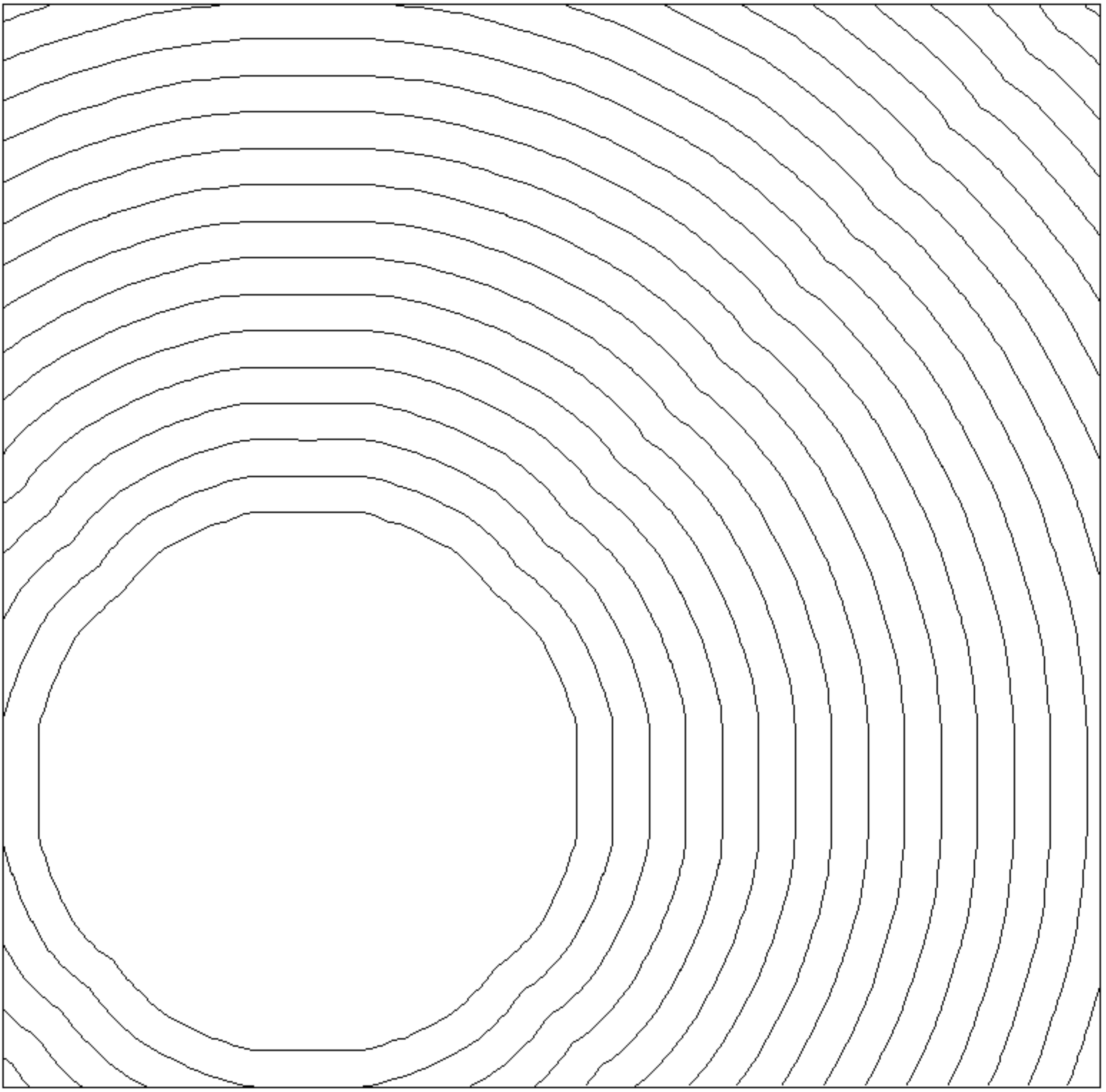} \\ 
(a)  &  (b) \\
\end{tabular}
\caption{Patchy solution for the dynamics Eikonal. (a) Value function and (b) a detail of its level sets. Small hollows are visible in correspondence of the lines $\{x=y\}$ and $\{x=-y\}$ cfr.\ Fig.\ \ref{fig:patchy_errorLinf}-(a).}
\label{fig:patchy_error_T}
\end{center}
\end{figure}

Fig.\ \ref{fig:patchy_error_T} shows the value function $U_P$ and its level sets for the Eikonal test. Here we see small perturbations where patches meet. 
It is interesting to note that they meet forming a hollow and not a discontinuity.


\subsection{Comparison of CPU times}\label{patchy_vs_dd}
In this section we compare the patchy algorithm with the classical domain decomposition algorithm in terms of CPU time. In case of parallel computation, we will always refer to the \textit{wall clock time}, and not to the sum of CPU times devoted to the task by each CPU running it.

Let us first explain why we expect the patchy algorithm overcomes the domain decomposition algorithm, considering again the Eikonal case with $R$=8, see Fig.\ \ref{fig:patchy_decomposition}-(a). 
If we visit the nodes in a single predefined order (i.e.\ we do \emph{not} implement the Fast Sweeping technique \cite{Z05} or similar ones), the eight subdomains need a different number of iterations to reach convergence. This is due to the fact that for some of them the visiting order corresponds to the upwind direction, while for the other subdomains the visiting order corresponds to the downwind direction. If we do not know \emph{a priori} that the eight subdomains are invariant with respect to the optimal dynamics, we cannot stop the computation in a subdomain before computations in \emph{all} subdomains are fully completed, because in any moment a new information can enter, making necessary new computations. On the contrary, if we know \emph{a priori} that subdomains do not depend on each other, we can safely stop the computation in a subdomain as soon as the solution reached convergence. Note that this argument is not related to the parallelization and holds also if only one core is employed. In addition, the higher the number of subdomains the more efficient is the computation. 

In Table \ref{tab:tempi_parziali} we report the execution times (in seconds) for the single steps of the patchy algorithm. 
Times for the Step 3.1 are the most interesting ones because this step is expected to be the slowest one after Step 3.2 (main computation on the fine grid). Thus, time spent in Step 3.1 could completely neutralize the advantage we hope to get in the subsequent main computation. 
As we can see, Step 3.1 is much more costly than Steps 1 and 2, but not so much compared with the main computation.
\begin{table}[h!]
\caption{CPU time. Cores: 4, $N_c$=32, Grid: $\Ncoarse=100^2\rightarrow N=800^2$, $R$=16}
\label{tab:tempi_parziali}
\begin{center}
\begin{tabular}{|c|c|c|c|c|}
\hline $ $ & Step 1 & Step 2 & Step 3.1 (all $j$'s)  & Step 3.2 (all $j$'s)
\\ \hline
\hline Eikonal & 2 & 1 & 23 & 409 \\
\hline Fan     & 2 & 2 & 52 & 796 \\
\hline Zermelo & 2 & 1 & 30 & 512 \\
\hline
\end{tabular}
\end{center}
\end{table}

Tables \ref{tab:comparison_eikonal16}, \ref{tab:comparison_eikonal32}, \ref{tab:comparison_fan} and \ref{tab:comparison_zermelo} report the CPU time (in seconds)  
for the three dynamics of Table \ref{tab:test2d} as a function of the number of cores (1,2,4) and the number of patches ($R$=2,4,8,16). They also report the speed-up obtained by 
the parallelization Method 1 (see Section \ref{sec:patchy}). For the Eikonal test we also vary the number of discrete controls ($N_c$=16 and 32). 
These results are compared with the best outcome of the domain decomposition method obtained varying the number of domains (again 2,4,8,16).\footnote{The CPU time of the domain decomposition method does not vary a lot varying the number of domains, but small differences are present. 
They are due to the different order in which nodes are visited and synchronization overhead at the end of each iteration.}
\begin{table}[h!]
\caption{CPU time (speed-up). Dynamics: Eikonal. $N_c$=16. Grid: $\Ncoarse=100^2\rightarrow N=800^2$}
\label{tab:comparison_eikonal16}
\begin{center}
\begin{tabular}{|l|l|l|l|l|l|}
\hline $ $ & $R=2$ & $R=4$ & $R=8$ & $R=16$ & Best DD
\\ \hline
\hline 1 core  & 1547 & 1076 & 1058 & 933 & 1571\\
\hline 2 cores & 845 (1.83) & 595 (1.81) & 574 (1.84)  & 504 (1.85) & 820 (1.92) \\
\hline 4 cores & 459 (3.37) & 325 (3.31) & 317 (3.34) & 271 (3.44) & 415 (3.79)\\
\hline
\end{tabular}
\end{center}
\end{table}
\begin{table}[h!]
\caption{CPU time (speed-up). Dynamics: Eikonal. $N_c$=32. Grid: $\Ncoarse=100^2\rightarrow N=800^2$}
\label{tab:comparison_eikonal32}
\begin{center}
\begin{tabular}{|l|l|l|l|l|l|}
\hline $ $ & $R=2$ & $R=4$ & $R=8$ & $R=16$ & Best DD
\\ \hline
\hline 1 core  & 2702 & 1897 & 1843 & 1623 & 2785\\
\hline 2 cores & 1462 (1.85) & 998 (1.90)  & 968 (1.90)  & 872 (1.86)  & 1430 (1.95)\\
\hline 4 cores & 771 (3.50) & 532 (3.57)  & 514 (3.59) & 435 (3.73) & 716 (3.89)\\
\hline
\end{tabular}
\end{center}
\end{table}
\begin{table}[h!]
\caption{CPU time (speed-up). Dynamics: Fan. $N_c$=32. Grid: $\Ncoarse=100^2\rightarrow N=800^2$}
\label{tab:comparison_fan}
\begin{center}
\begin{tabular}{|l|l|l|l|l|l|}
\hline $ $ & $R=2$ & $R=4$ & $R=8$ & $R=16$ & Best DD
\\ \hline
\hline 1 core  & 3712 & 3322 & 3049 & 3172 & 4163\\
\hline 2 cores & 2020 (1.84)& 1746 (1.90)& 1596 (1.91)& 1559 (2.03)& 2124 (1.96)\\
\hline 4 cores & 1032 (3.60)& 900  (3.69)& 841 (3.63)& 852 (3.72)& 1069 (3.89)\\
\hline
\end{tabular}
\end{center}
\end{table}
\begin{table}[h!]
\caption{CPU time (speed-up). Dynamics: Zermelo. $N_c$=32. Grid: $\Ncoarse=100^2\rightarrow N=800^2$}
\label{tab:comparison_zermelo}
\begin{center}
\begin{tabular}{|l|l|l|l|l|l|}
\hline $ $ & $R=2$ & $R=4$ & $R=8$ & $R=16$ & Best DD
\\ \hline
\hline 1 core  & 3113 & 2675 & 2126 & 2018 & 3209\\
\hline 2 cores & 1651 (1.89)& 1404 (1.91)& 1111 (1.91)& 1054 (1.91)& 1640 (1.96)\\
\hline 4 cores & 871  (3.57)& 721  (3.71)& 584  (3.64)& 545  (3.70)& 825 (3.89)\\
\hline
\end{tabular}
\end{center}
\end{table}

We see that the speed-up is very satisfactory and proves that the parallelization Method 1 we implement here is sound. Moreover, we see that the CPU time decreases remarkably as the number of patches $R$ increases. 
For $R$=16 the CPU time is considerably smaller than that of the best domain decomposition method. This is one of the main results of the paper. 

Differences among $N_c$=16 and $N_c$=32 are instead less clear, although the patchy algorithm should have an advantage for large $N_c$ because of the smaller ratio between CPU time 
for Step 3.1 (one discrete control) and Step 3.2 ($N_c$ discrete controls).

\medskip
\textit{Remark 4.1.} The Fast Sweeping technique can mitigate the performances of patchy method, since it clears the differences between domains with the same number of nodes, 
but it cannot neutralize them completely. Indeed, the patchy algorithm has the clear advantage that no synchronization or crossing information among processors are needed. 
This is a great advantage when using distributed-memory parallel computers (for which Method 2 is designed), where communications are performed via cables connecting cluster nodes. 
This advantage is not really included in our experiments because our cores share a common RAM.

\subsection{Patchy method with obstacles} 
\begin{figure}[h!]
\begin{center}
\begin{tabular}{cc}
\includegraphics[width=0.43\textwidth]{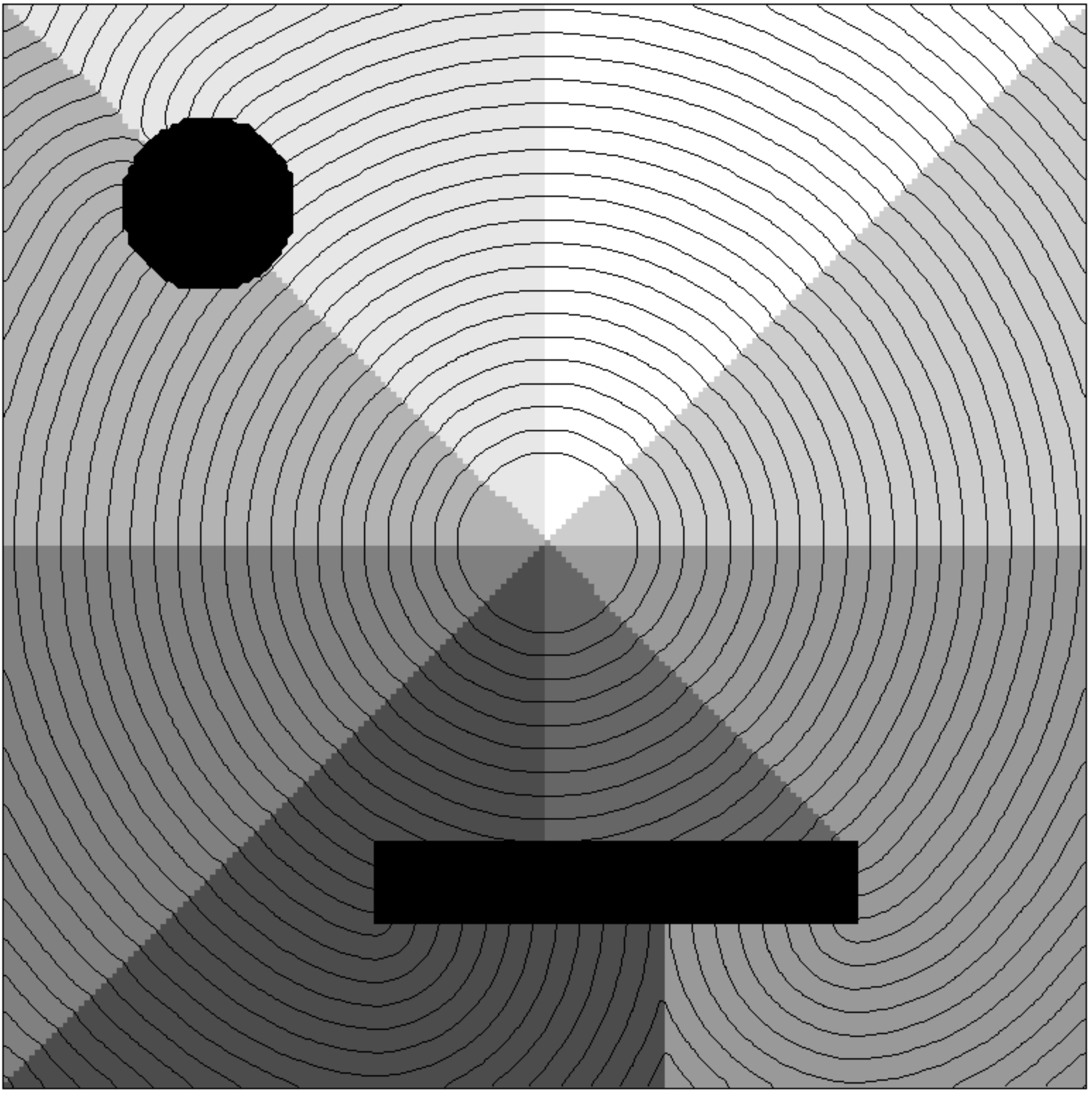} &
\includegraphics[width=0.50\textwidth]{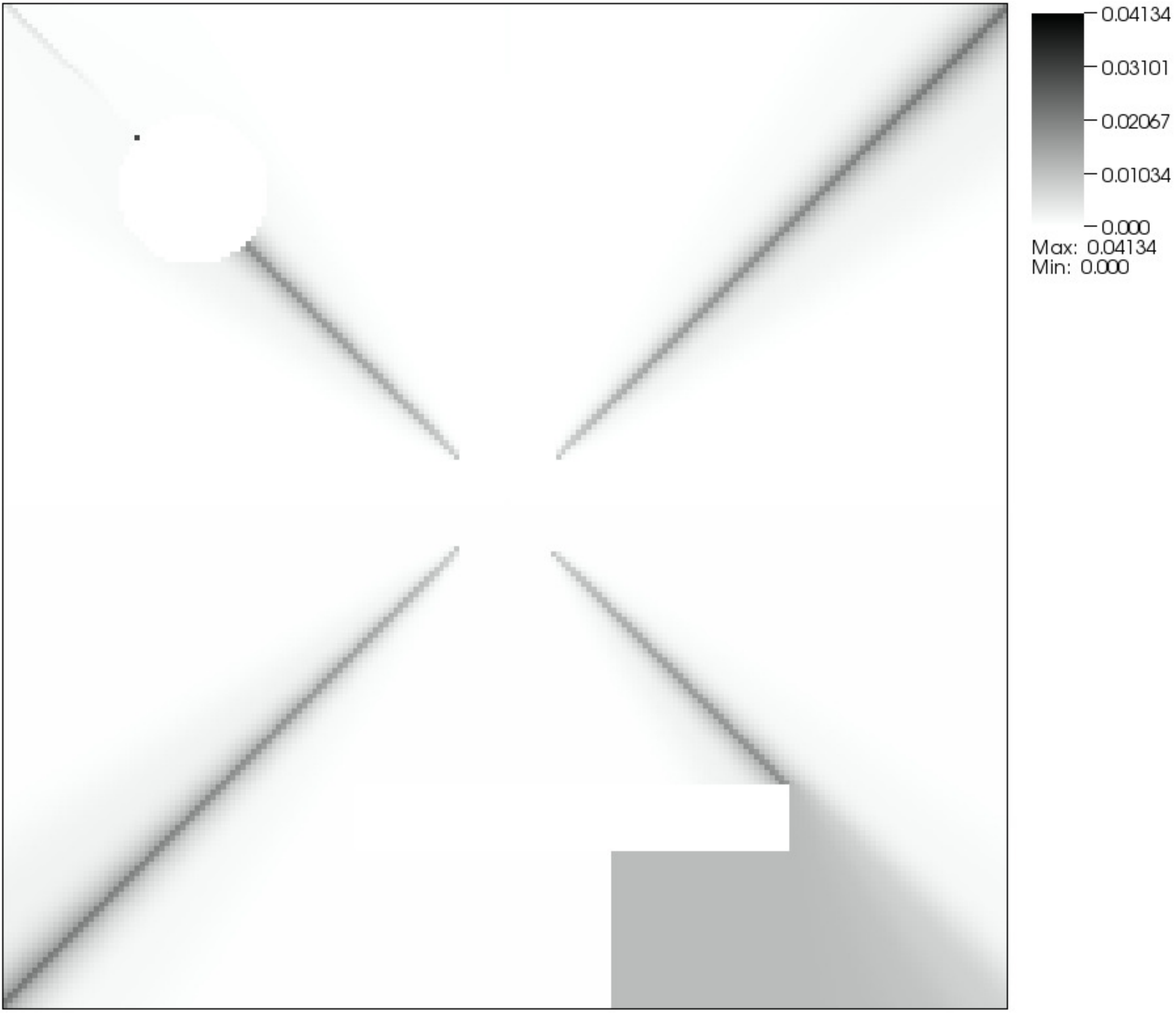} \\ 
(a)  & \hskip-22pt (b) 
\end{tabular}
\caption{Patches bypass the obstacles driven by the dynamics. (a) Obstacles (black), level sets of the value function and patchy decomposition. (b) Patchy error $E$.}
\label{fig:patchy_obstacles}
\end{center}
\end{figure}

We have also tried to use the patchy algorithm to solve a minimum time problem with Eikonal dynamics and obstacles. In Fig.\ \ref{fig:patchy_obstacles} we show the obstacles (one circle and one rectangle), the level sets of the solution, the patchy decomposition and the patchy error $E$. The behaviour of the patchy decomposition is correct because the dynamics drives the patches around the obstacles. If not influenced by the obstacles, the error is concentrated around the boundaries of the patches as expected. 
Instead, when a boundary meets an obstacles, the error can either stop propagating (see the circle) or spread out (see the rectangle, right side).

\subsection{Limitations of the patchy method}\label{limitations}
The overall efficiency of the patchy method depends on the dynamics and on the shape of the target $\Omega_0$. Unfortunately, it is not always possible 
to get a suitable patchy decomposition which allows one to run the algorithm. This can happen for example if the target is very small and then it cannot be divided in $R$ subdomains. 

Another issue, much more difficult to fix, comes out whenever there is a large difference between the sizes of the patches and possibly some of them degenerate in a subset of a few grid nodes. This is the case of the classical ``Lunar Landing'' problem
$$
d=2\,,\qquad f(x_1,x_2,a)=(x_2,a)\,,\qquad A=\{-1,1\}\,, \qquad \Omega_0=B_2(0,\varepsilon).
$$
In this case the patchy decomposition consists of 2 large domains and $R-2$ smaller domains, see Fig.\ \ref{fig:krener}-(a). 
The small domains degenerate to sets of dimension one when $\varepsilon$ tends to zero, because all the optimal trajectories tend to meet in only two switching lines.

A third dangerous case arises when some regions in $\Omega_0$ are not reachable. If the dynamics make it impossible to reach the target from some point $x\in\Omega$, 
the value function is set to $u(x)=+\infty$. From the numerical point of view, the solution stays frozen at the value given as initial guess. 
At these points the optimal control (\ref{feedbackcontrol}) is not uniquely defined, and then the patchy decomposition cannot be build. 
On the other hand, in the not-reachable regions the solution is in some sense already computed, then the issue can be easily fixed by a slight modification of the algorithm. 
After Step 1 we locate the regions where ${\widetilde U}_P$ is very large and then we do not consider those regions in the rest of computations.
\subsection{Patchy decomposition for non-target problems}\label{krener} 
\begin{figure}[b!]
\begin{center}
\begin{tabular}{cc}
\includegraphics[width=0.43\textwidth]{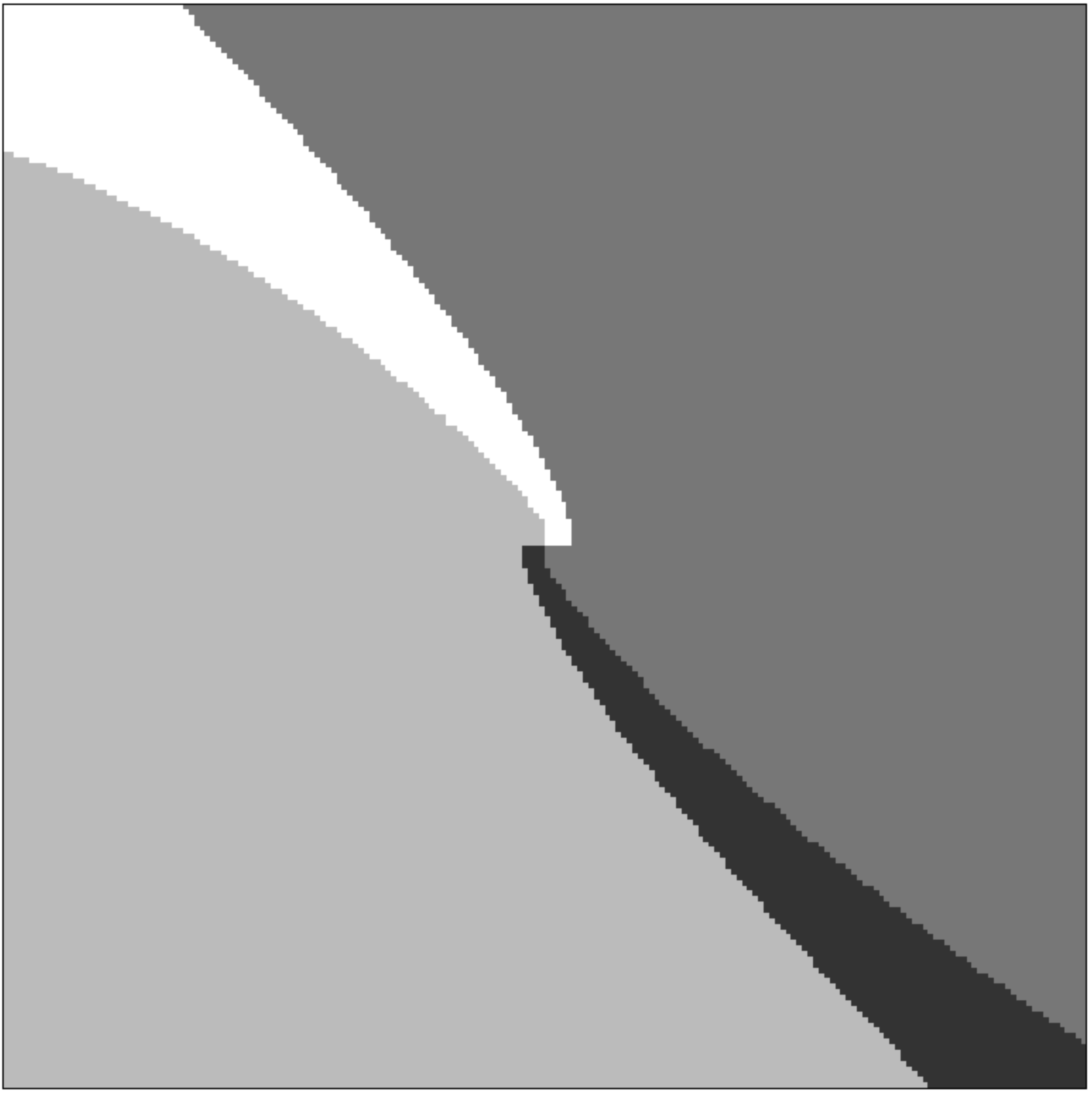} &
\includegraphics[width=0.50\textwidth]{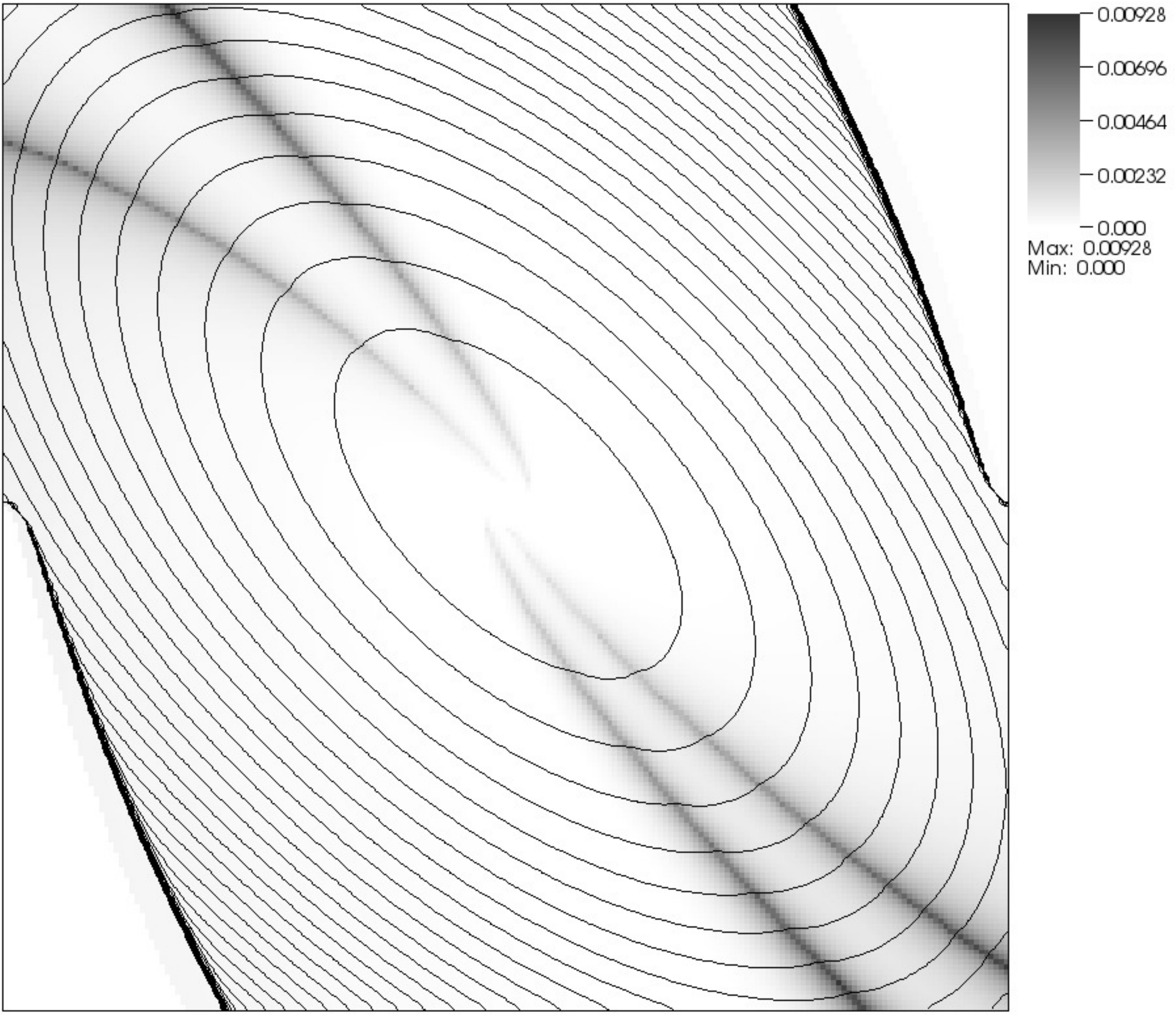} \\
(a)  &  \hspace{-10 mm} (b)
\end{tabular}
\begin{tabular}{cc}
\hspace{-10 mm} \includegraphics[width=0.43\textwidth]{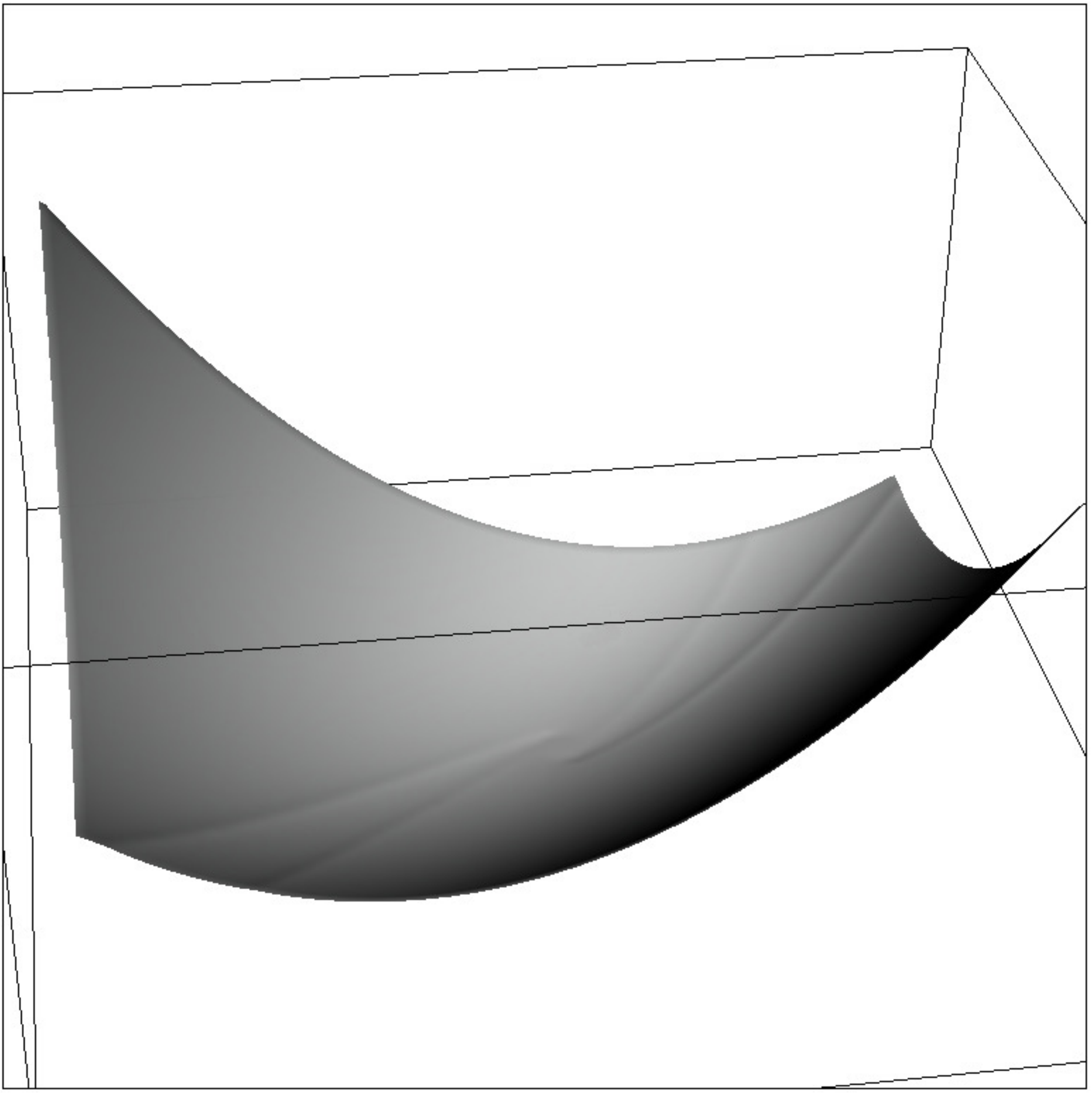} &
\hspace{-0 mm} \includegraphics[width=0.43\textwidth]{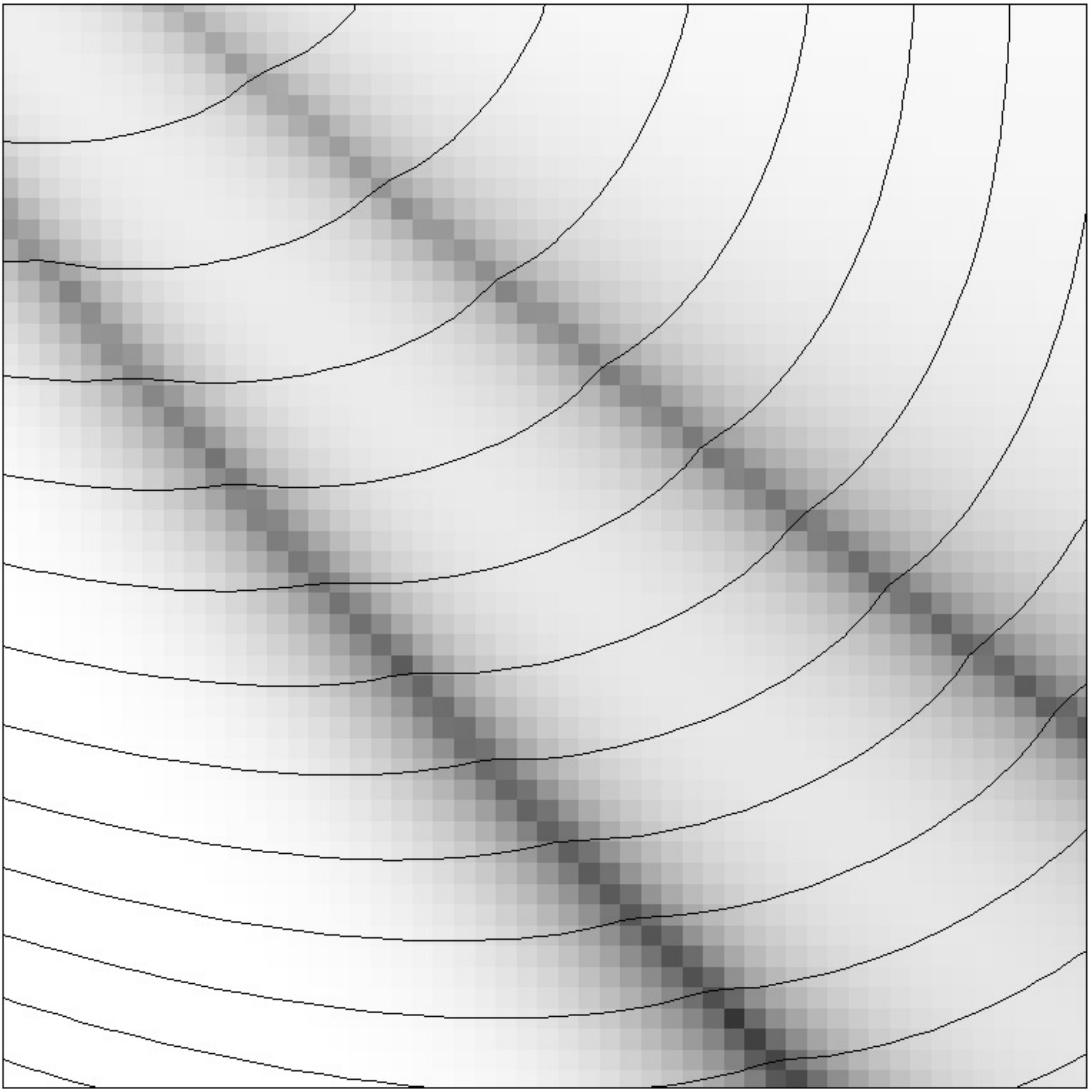} \\
(c)  &  \hspace{-0 mm} (d)
\end{tabular}
\caption{
(a) Decomposition in $R$=4 patches, 
(b) error function and level sets of the patchy solution (the solution is truncated at value $u=3$ in order to remove the boundary effects which are very important for this dynamics), 
(c) patchy solution $U_P$, (d) a detail of the error function shown in (b).}
\label{fig:krener}
\end{center}
\end{figure}
In the case of non-target problems we cannot in principle build the patchy decomposition. Indeed, we recall that our patchy decomposition starts from a decomposition of the target $\Omega_0$. 
Nevertheless, in some special situations it is still possible to achieve the patchy decomposition, using some a-priori knowledge on the solution of the problem. 
This is the case of the infinite horizon problem associated to the \emph{linear-quadratic regulator}, studied by Krener and Navasca in \cite{KN08} as a test for their patchy method 
(see also \cite{KN07} and the Introduction for a brief description):   
\begin{equation}\label{LQR}
\min_{a\in\R}\int_0^\infty \left(\frac12 (y_1^2+y_2^2)+\frac12 a^2 \right)dt \quad\textrm{ subject to } \quad
\left\{\begin{array}{l}
\dot y_1=y_2 \\
\dot y_2=a
\end{array}
\right. 
\end{equation}
with $y_1(t=0)=x_1$ and $y_2(t=0)=x_2$. The exact value function for this problem is
$$
u(x_1,x_2)=\frac12
\left(\!\!\begin{array}{c} x_1 \\ x_2 \end{array}\!\!\right)
\left(\!\!\begin{array}{cc} \sqrt{3} & 1 \\ 1 & \sqrt{3} \end{array}\!\!\right)
(x_1 \ x_2) 
$$
and it is easy to check that the origin $(0,0)$ is the only source of all characteristic curves. 
Then, using a small ball $B_2(0,\varepsilon)$ as a fictitious target, we are able to generate the patchy decomposition and then run the algorithm normally. 
In Fig.\ \ref{fig:krener} we show the outcome of the simulation with the following parameters:
$\Omega=[-1,1]^2$, $\varepsilon=0.05$, $R$=4, $N_c$=101, $A=[-3,3]$, $\Ncoarse$=$100^2$, $N$=$200^2$. Note that the choice $R$=4 is due to limitations discussed in Section \ref{limitations}.
Moreover, we cannot impose state constraint boundary conditions at the boundaries of the patches, since at some nodes the dynamics in (\ref{LQR}) does not point inside the patch they belong to for any $a\in A$. As already discussed in Section \ref{sec:patchy}, this issue is solved by using $U_P^{(0)}$ as Dirichlet boundary condition:
\begin{equation}\label{dbc}
U_P|_{\partial \Omega^j}= U_P^{(0)}|_{\partial \Omega^j}\,,\qquad j=1,\dots,R\,.
\end{equation}
Patchy errors are $E_1$=0.00012 and $E_\infty$=0.009. Note that they are generally smaller than those in Tables \ref{tab:patchy_error_eikonal}-\ref{tab:patchy_error_zermelo} 
(computed for other dynamics) because of the more favorable boundary condition \eqref{dbc}.


\section{Patchy method's add-ons}
The patchy algorithm proposed in Section \ref{sec:patchy} has a multigrid nature, meaning that the computation of the solution on a rough grid is needed to start the optimal domain decomposition. 
Once this preliminary effort is done, it appears to be natural to use all the information we have collected in order to speed up the algorithm. 
First of all, in the next tests we impose by default the boundary condition \eqref{dbc} (note that it becomes available only after the computation on the rough grid). Further multigrid advantages we can take into account are listed in the following. 
\begin{enumerate}
\item[AO1.] We use $U_P^{(0)}$ computed in Step 2 as initial guess for Step 3.2. In this way we save some iterations to reach convergence.\\
\item[AO2.] Before Step 3.2 we order the nodes belonging to each patch in such a way they fit as much as possible the \textit{causality principle} \cite{SV03}. 
For example we can order the nodes with respect to their values. This ordering is optimal if the characteristic lines coincide with the gradient lines of the solution, as it happens in the case of the Eikonal equation. 
In general this is not true, anyway this ordering is often not too far from the optimal one. \\
\item[AO3.] In Step 3.2 we reduce the number of discrete controls used in the numerical scheme, 
eliminating those controls which are ``far'' from the optimal one $a^*_{\widetilde k}(x_i)$ as computed by the first computation on the rough grid (Step 2). 
For example, if $A=B_2(0,1)$ we can introduce a reduction factor $r>1$ and replace $A$ with the set 
$$
A_{r}=\left\{a\in A\ :\ a\cdot a^*_{\widetilde k}\geq \cos\left(\frac{\pi}{r}\right)\right\}.
$$
This is the only add-on which introduces a new error in the solution, anyway it is negligible in most cases.\\
\end{enumerate}

We point out that the patchy method can easily become an actual multigrid method. Indeed, we can in principle repeat the algorithm introducing a sequence of grids $G_1,G_2,\ldots$ one finer than the other, until the desired precision is reached.

In order to study the effect of the previously described add-ons, we introduce them separately and we compare the CPU time with the basic algorithm. Then, we apply all the features together. Results are reported in Table \ref{tab:addons}.
\begin{table}[h!]
\caption{Effects of add-ons. Cores: 2, $R=8$, $N_c=32$. Controls reduced by factor 4}
\label{tab:addons}
\begin{center}
\begin{tabular}{|c|c|c|c|c|c|c|}
\hline dynamics & grid size & basic & AO1 & AO2 & AO3 & AO1+AO2+AO3
\\ \hline
\hline Eikonal & $100^2\rightarrow 200^2$ & 20.0   & 19.2   & 9.6   & 9.1   & 5.7   \\
\hline Eikonal & $100^2\rightarrow 400^2$ & 130.7  & 130.2  & 40.5  & 43.6  & 17.8  \\
\hline Eikonal & $100^2\rightarrow 800^2$ & 928.1  & 924.6  & 238.8 & 298.1 & 100.6 \\
\hline Fan     & $100^2\rightarrow 200^2$ & 31.9   & 31.0   & 11.4  & 14.0  & 7.6   \\
\hline Fan     & $100^2\rightarrow 400^2$ & 209.8  & 205.7  & 43.5  & 72.3  & 20.6  \\
\hline Fan     & $100^2\rightarrow 800^2$ & 1571.9 & 1564.0 & 247.3 & 529.6 & 110.6 \\
\hline Zermelo & $100^2\rightarrow 200^2$ & 23.2   & 22.6   & 11.5  & 10.7  & 6.7   \\
\hline Zermelo & $100^2\rightarrow 400^2$ & 143.5  & 142.4  & 46.2  & 51.0  & 20.3  \\
\hline Zermelo & $100^2\rightarrow 800^2$ & 1071.4 & 1057.9 & 290.1 & 345.5 & 111.3 \\
\hline
\end{tabular}
\end{center}
\end{table}

Note that CPU times for this test are lower than those in Section \ref{patchy_vs_dd} because of the more favorable boundary condition \eqref{dbc}.

\section{Numerical tests in dimension three}\label{testsHD}
We solve the three 3D minimal time problems of the form \eqref{HJB} listed in Table \ref{tab:test3d}. The numerical domain is $\Omega=[-2,2]^3$ for all tests.
\begin{table}[h!]
\caption{Three-dimensional numerical tests}
\label{tab:test3d}
\begin{center}
\begin{tabular}{|l|l|l|l|l|}
\hline Name       & $d$ & $f(x_1,x_2,x_3,a)$   & $A$ & $\Omega_0$ \\ \hline
\hline Eikonal 3D & 3   & $a$                  & $B_3(0,1)$ & $B_3(0,0.5)$ \\
\hline Fan 3D     & 3   & $|x_1+x_2+x_3+0.1|a$ & $B_3(0,1)$ & $\{x_1=0\}$ \\
\hline Brockett 3D \cite{BRM92, MB99} & 3   & $(a_1,a_2,x_1a_2-x_2a_1)$  & $[-5,5]^2$ & $B_3(0,0.25)$ \\
\hline
\end{tabular}
\end{center}
\end{table}
For the first two tests we used $N_c$=189 discrete controls uniformly distributed on the unit sphere\footnote{In order to define an uniform distribution of discrete controls on the unit sphere, we used the vertices of a geosphere obtained by recursion starting from an icosahedron.}, 
while for the last test we used only $N_c$=9 discrete controls in $\{-5,0,5\}^2$. 
The latter choice is motivated by the fact that using a larger number of discrete controls in $[-5,5]^2$ does not lead to a different result, since the optimal strategy always requires to saturate the control to the extremal admissible values ($\pm 5$, in this case).

\begin{figure}[h!]
\begin{center}
\includegraphics[width=0.50\textwidth]{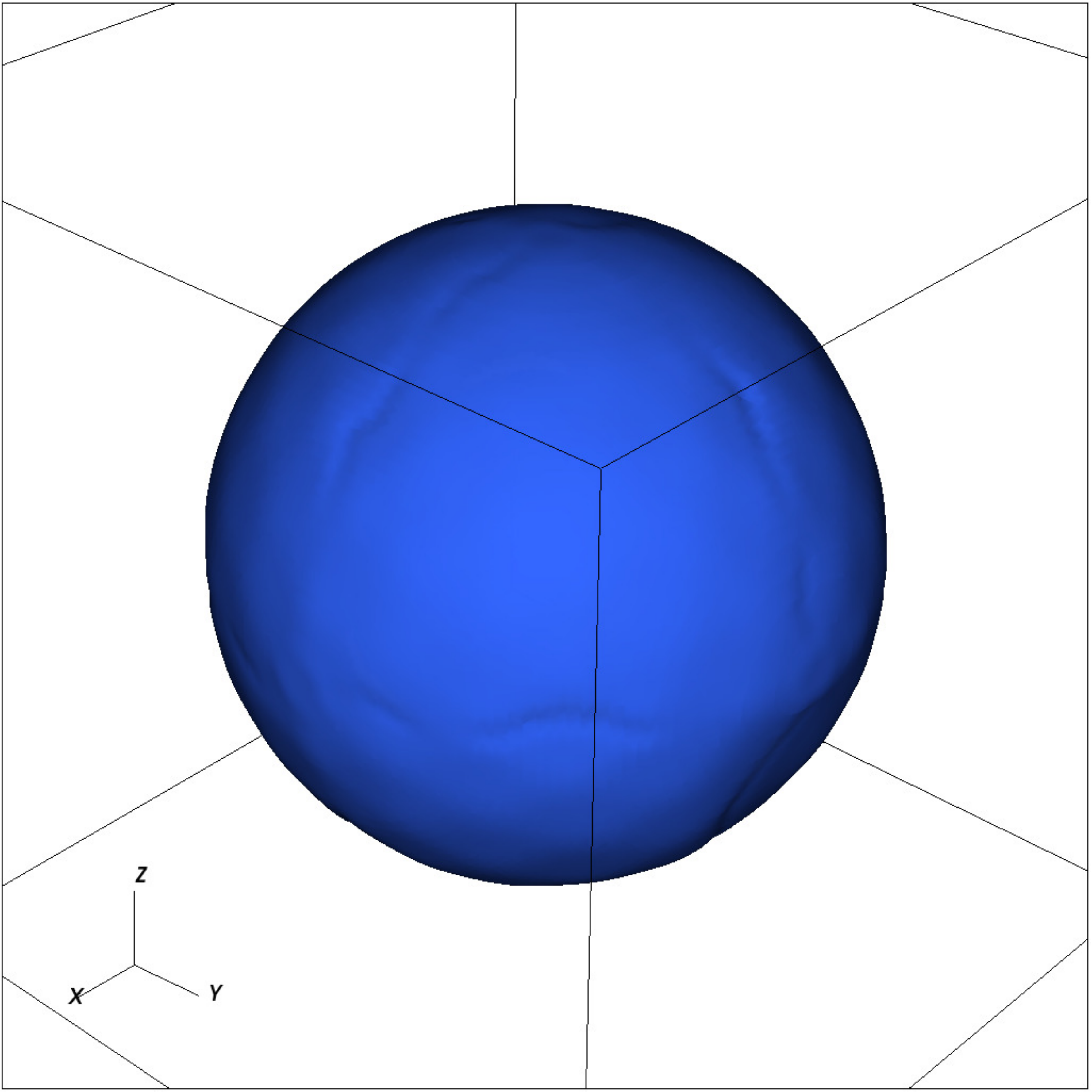}
\caption{One level set of the value function for the Eikonal 3D dynamics}
\label{fig:3Deikonal}
\end{center}
\end{figure}

\begin{figure}[h!]
\begin{center}
\begin{tabular}{cc}
\includegraphics[width=0.45\textwidth]{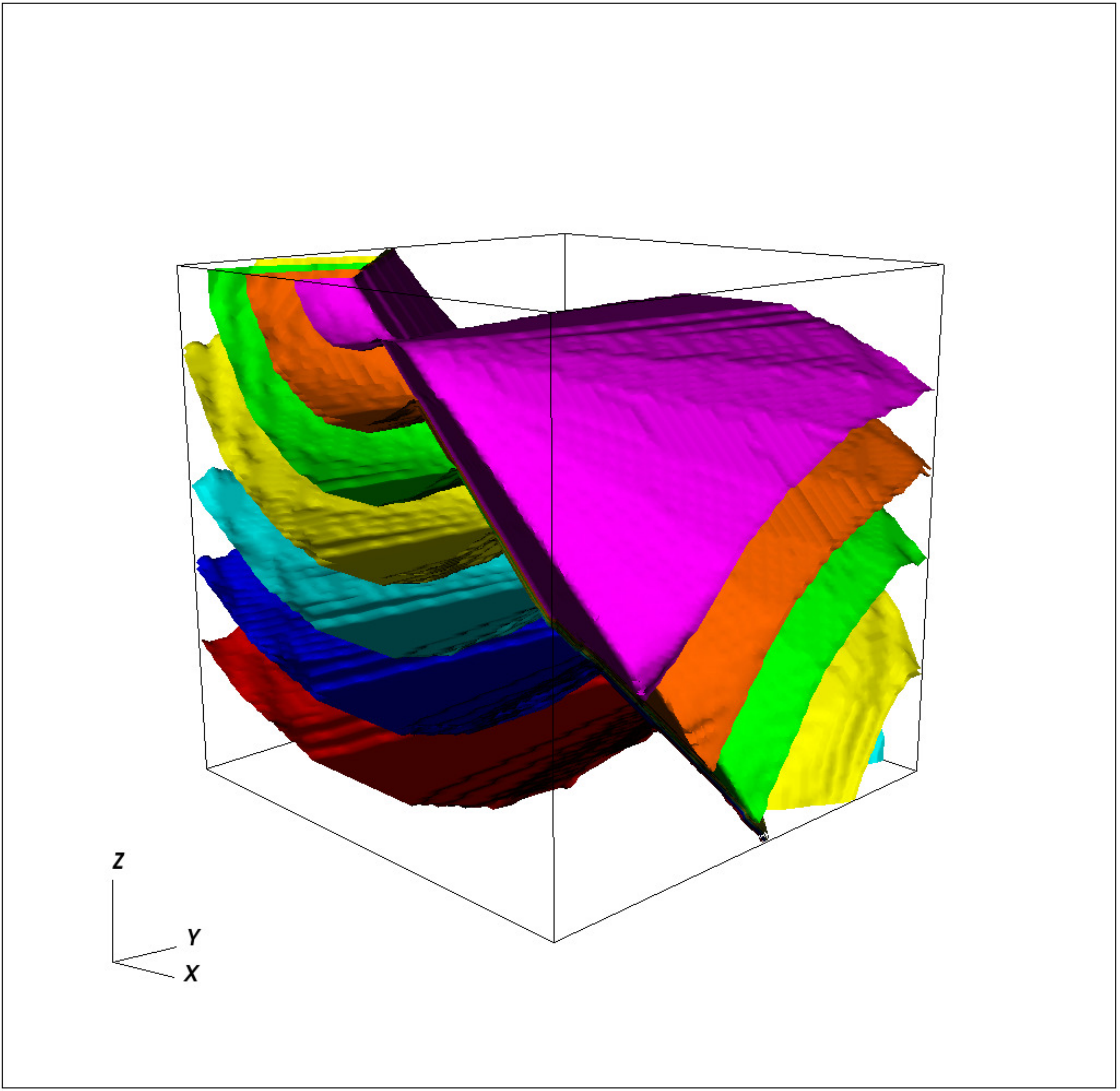} &
\includegraphics[width=0.45\textwidth]{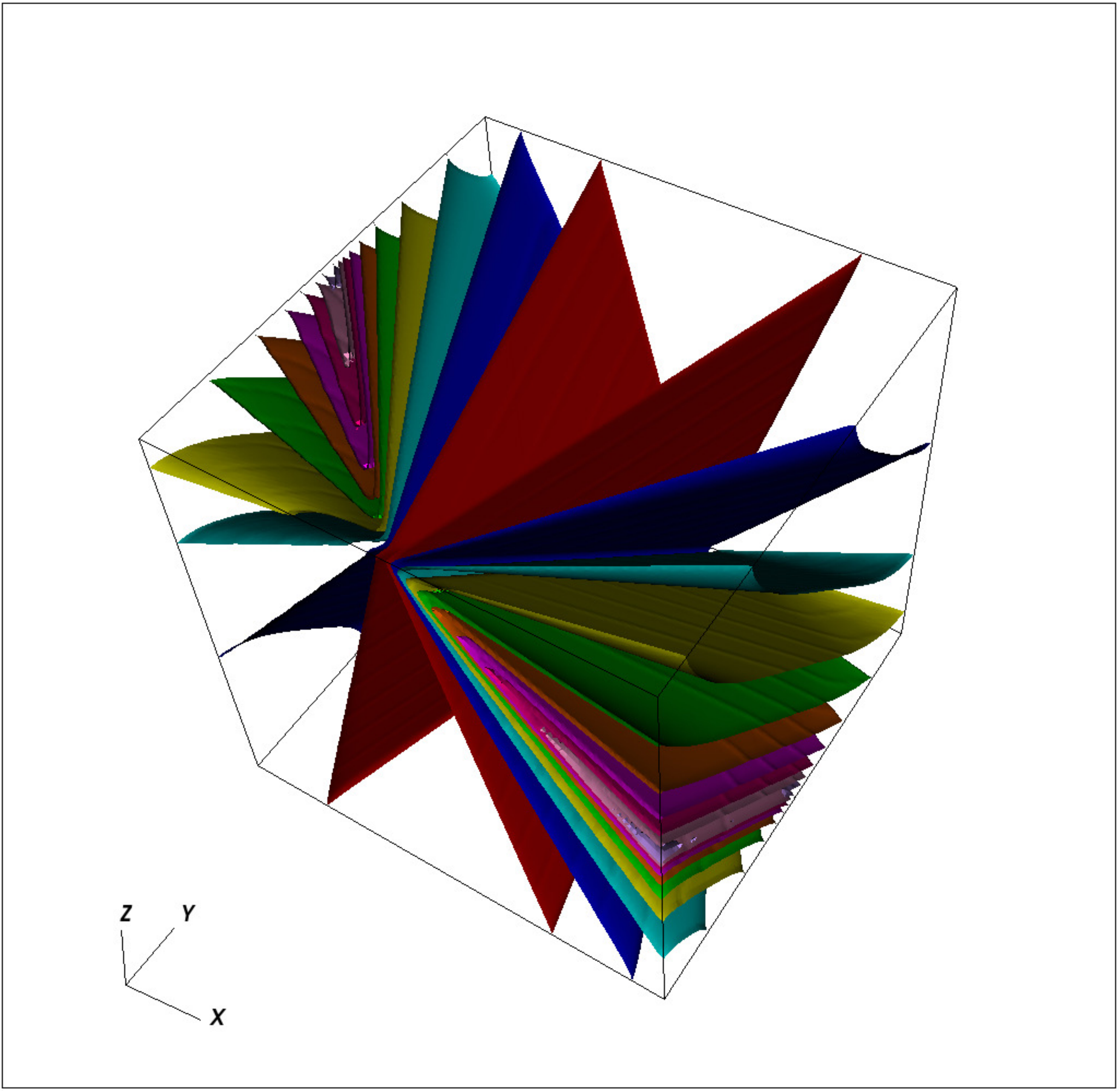} \\
(a)  &  (b)
\end{tabular}
\begin{tabular}{c}
\includegraphics[width=0.50\textwidth]{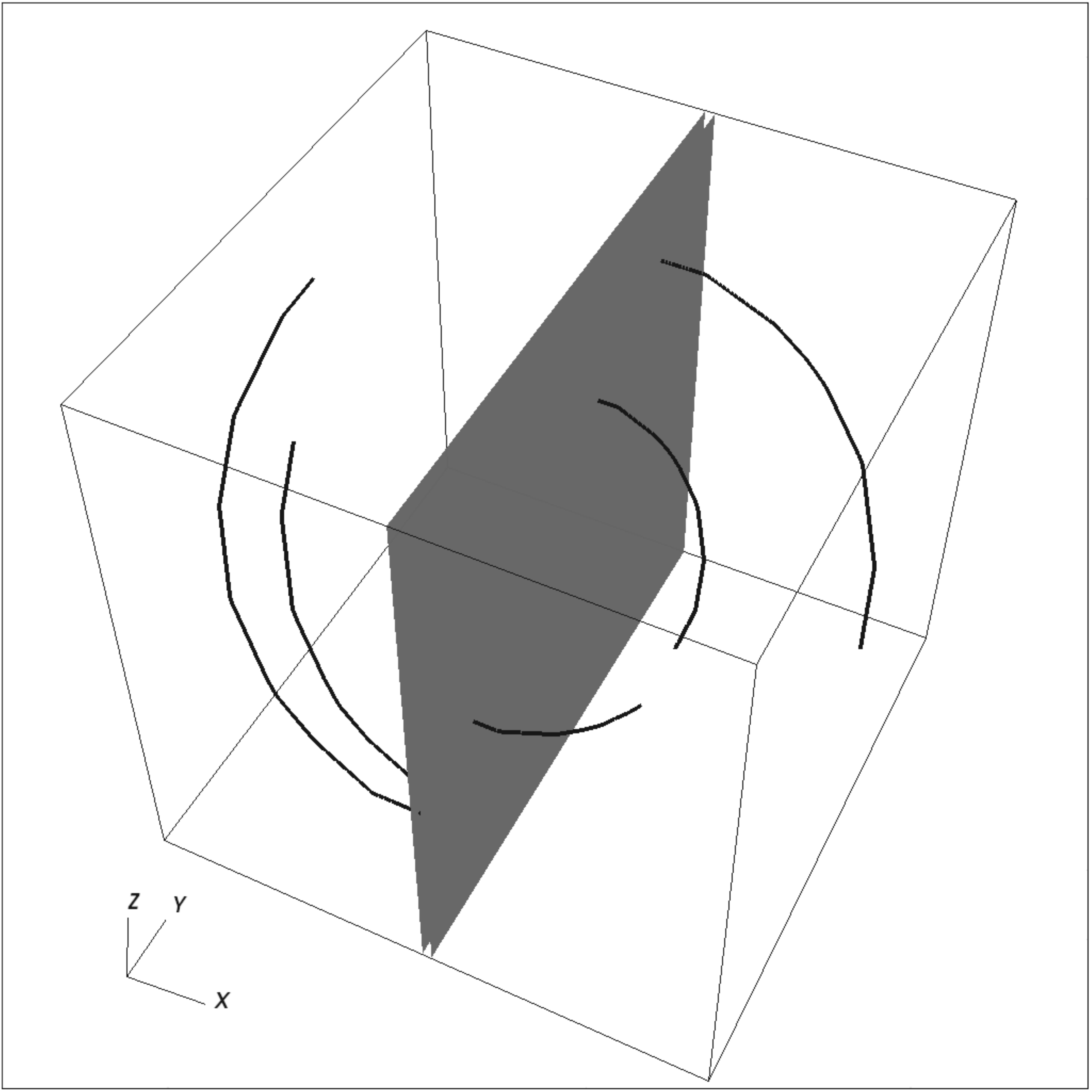} \\
(c)  \\
\end{tabular}
\caption{Fan 3D: (a) Decomposition with 8 patches, (b) level sets of the solution, and (c) some optimal trajectories to the target}
\label{fig:3Dfan}
\end{center}
\end{figure}

Results are reported in Table \ref{tab:3D}. 
\begin{table}[h!]
\caption{3D tests. Cores: 4, $R=8$. Add-ons enabled (Eikonal and Fan: controls reduced by factor 4, Brockett: not reduced)}
\label{tab:3D}
\begin{center}
\begin{tabular}{|c|c|c|c|c|}
\hline dynamics & grid size & CPU time & $E_1$ & $E_\infty$
\\ \hline
\hline Eikonal 3D   & $50^3\rightarrow 100^3$ & 183   & 0.00052  & 0.035  \\
\hline Eikonal 3D   & $50^3\rightarrow 200^3$ & 1217  & 0.00045  & 0.042  \\
\hline
\hline Fan 3D       & $50^3\rightarrow 100^3$ & 165   & 0.00100  & 0.187  \\
\hline Fan 3D       & $50^3\rightarrow 200^3$ & 1269  & 0.00087  & 0.305  \\
\hline
\hline Brockett 3D  & $50^3\rightarrow 100^3$ & 132   & 0.00358  & 0.024  \\
\hline Brockett 3D  & $50^3\rightarrow 200^3$ & 1557  & 0.00258 & 0.020  \\
\hline
\end{tabular}
\end{center}
\end{table}
Considering the large number of discrete controls used for Eikonal 3D and Fan 3D, the CPU time is remarkable. 
Fig.\ \ref{fig:3Deikonal} shows a level set of the value function for the Eikonal 3D dynamics. It is perfectly visible that error is located where the patches meet.
Fig.\ \ref{fig:3Dfan} shows instead the results for Fan 3D with $200^3$ nodes. In Figs.\ \ref{fig:3Dfan}(a,b) we show the boundaries of the patches and some level sets of the solution, respectively. 
Level sets should be plans, but the state constraints imposed by the computational box $\Omega$ bend them near $\partial\Omega$. 
In Fig.\ \ref{fig:3Dfan}(c) we show some optimal trajectories to the target.

\begin{figure}[t!]
\begin{center}
\begin{tabular}{cc}
\includegraphics[width=0.45\textwidth]{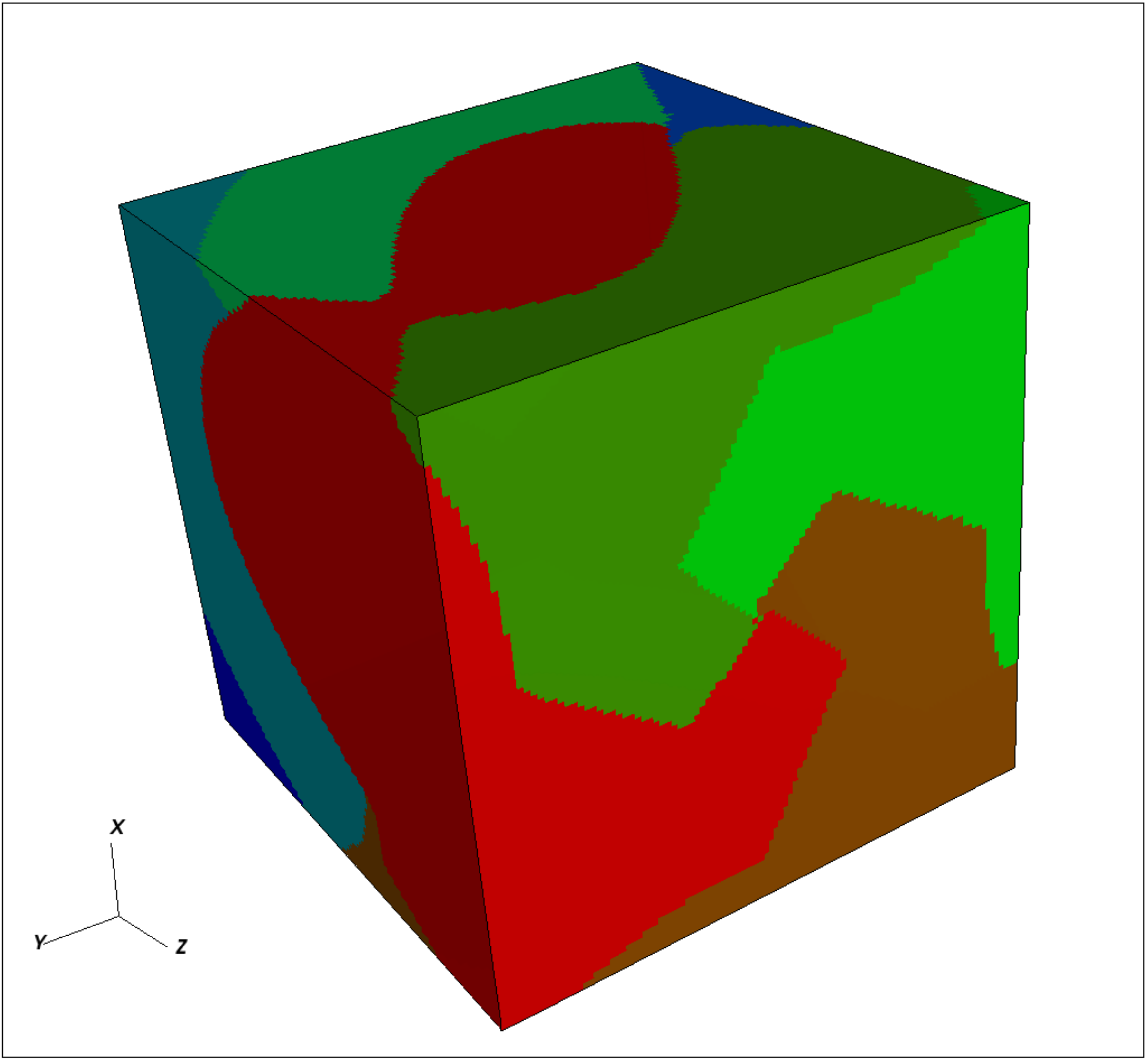} &
\includegraphics[width=0.45\textwidth]{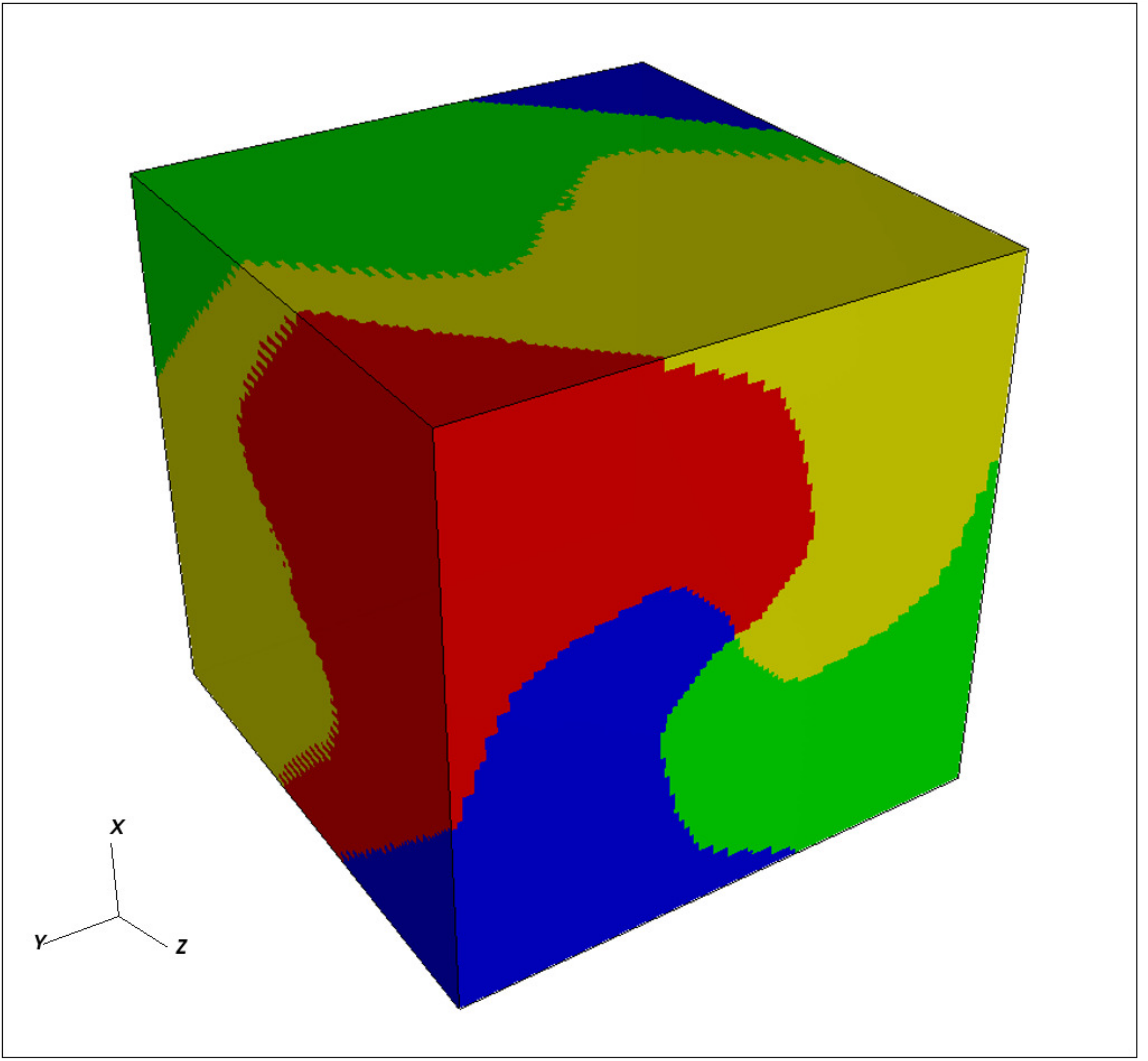} \\
(a)  &  (b)
\end{tabular}
\begin{tabular}{c}
\includegraphics[width=0.50\textwidth]{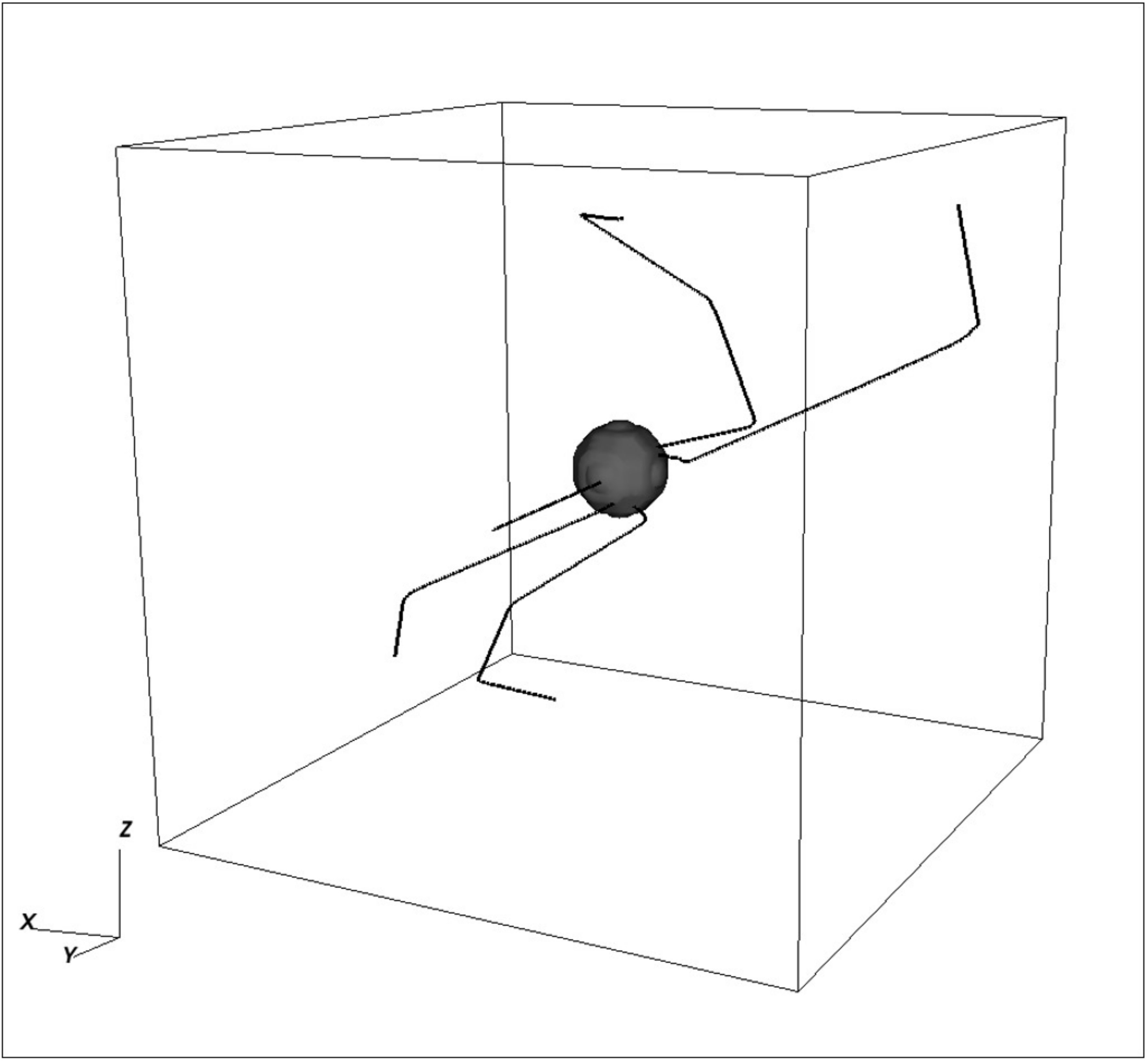} \\
(c)  \\
\end{tabular}
\caption{Brockett 3D: (a) Decomposition with 8 patches, (b) Regions with constant optimal controls, and (c) some optimal trajectories to the target}
\label{fig:3Dbrockett}
\end{center}
\end{figure}

Results for Brockett problem are different from the previous tests. First, CPU time turns out to be high with respect to the small number of discrete controls in use (just $N_c=9$ controls).
This could be related to the fact that characteristics are broken lines (see Fig.\ \ref{fig:3Dbrockett}(c)) that do not go directly to the target as in the Eikonal equation, 
nor bend slightly as for the Fan dynamics (see Fig.\ \ref{fig:3Dfan}(c)), but change direction instantaneously (see also control switch regions in Fig.\ \ref{fig:3Dbrockett}(b)), so that this dynamics takes much more time to move information through the domain. 
Second, the patchy error $E_1$ is quite large if compared to the other dynamics (see Table \ref{tab:3D}). 
This depends on the fact that the patchy decomposition obtained for this dynamics is rather complicated (see Fig.\ \ref{fig:3Dbrockett}(a)), in particular patches arrange themselves in (suggestive) sets with very large boundary areas, and this increases the number of nodes with large error $E$. 
\section*{Concluding remarks and future directions}
In this paper we have proposed a new numerical method for optimal control problems which tries to mimic the ``patchy decomposition'' proposed by Ancona and Bressan \cite{AB99}. We have investigated the serial implementation of the algorithm as well as one of the possible ways to parallelize it (Method 1), particularly suitable for shared-memory architectures. The new method is shown to be faster than the classical domain decomposition algorithm, since it avoids useless computations at nodes that have already reached convergence. At present, the main drawback of our approach is the fact that we have almost no control on the size of the patches, which only depends on the initial partition of the target and the dynamics.

Many points need to be investigated in the next future. The first is the parallelization on distributed-memory architectures, where patches are processed in parallel (Method 2). Moreover, we plan to improve the dynamics-dependent decomposition in such a way that the size of the patches is controllable. Finally, we aim at obtaining a convergence result for the scheme and, possibly, an a-priori estimate for the patchy error. These results, coupled with the previous results by Ancona and Bressan, will produce a constructive provably convergent method for some classes of optimal control problems.



\end{document}